\documentclass[12pt,a4paper,oneside,reqno]{amsproc}
\usepackage{amsmath} 
\usepackage{amssymb,amsfonts,mathrsfs}
\usepackage[colorlinks=true,linkcolor=blue,citecolor=blue]{hyperref}
\usepackage[driver=pdftex,lmargin=2.5cm,rmargin=2.5cm,heightrounded=true,centering]{geometry}
\usepackage{color}
\usepackage{epsfig}  
\usepackage{graphicx}  
\usepackage{xy}
\xyoption{all}
\begin{document}
\title[Fractional Prescribed Webster Scalar
Curvature Problem]
{Existence results for the fractional Q-curvature problem on three dimensional CR sphere}

\author{Chungen Liu}

\author{Yafang Wang}
\address{Nankai University, Tianjin 300071, P. R. China.}
\email{liucg@nankai.edu.cn(for Chungen Liu), wangyafang4727052@126.com(for Yafang Wang)}
\thanks{ The first author is partially supported by the NSF of China (11471170, 10621101), 973 Program of
MOST (2011CB808002) and SRFDP}
\newtheorem{theorem}{Theorem}[section]
\newtheorem{definition}{Definition}[section]
\newtheorem{lemma}{Lemma}[section]
\newtheorem{proposition}{Proposition}[section]
\newtheorem{corollary}{Corollary}[section]
\newtheorem{remark}{Remark}[section]
\renewcommand{\theequation}{\thesection.\arabic{equation}}
\catcode`@=11 \@addtoreset{equation}{section} \catcode`@=12
\newcommand{\optional}[1]{\relax}

\setcounter{secnumdepth}{3}
\setcounter{section}{0} \setcounter{equation}{0}
\numberwithin{equation}{section}
\newcommand{\MLversion}{1.1}
\begin{abstract}
In this paper the fractional Q-curvature problem on three dimensional  CR sphere is considered. By  using the  critical points theory at infinity, an existence result is obtained.
\end{abstract}
\maketitle

\def\Pf{\noindent{\it Proof.} }
\def\pf{\noindent{\it Proof.} }
\def\endd{\hfill $\Box$}
\def\lb{\label}
\def\re{\color{red}}
\def\bl{\color{blue}}
\def\ye{\color{yellow}}

\section{Introduction and main result}

The sphere $\mathbb{S}^{2n+1}$ is the boundary of the unit ball of $\mathbb{C}^{n+1}$. It is a contact manifold with a standard contact form $\theta_1$.  We denote by $(\mathbb{S}^{2n+1}, \theta)$ the contact sphere with its contact form $\theta$. Let $K: \mathbb{S}^{2n+1}\rightarrow \mathbb{R}$ be a $C^2$ positive function. The prescribed Webster scalar curvature problem
on $\mathbb{S}^{2n+1}$ is to find suitable conditions on $K$ such that $K$ is the Webster scalar curvature
for some contact form $\tilde{\theta}$ on $\mathbb{S}^{2n+1}$, CR equivalent to $\theta_1$. If we set $\tilde{\theta}=u^{\frac{2}{n}}\theta_1$ , where $u$ is a
smooth positive function on $\mathbb{S}^{2n+1}$, then the above problem is equivalent to solving the
following PDE:
\begin{align*}
\begin{cases}
L_{\theta_1}u=Ku^{1+\frac{2}{n}},\\
u>0,
\end{cases}
\end{align*}
where $L_{\theta_1}$ is the conformal Laplacian of $\mathbb{S}^{2n+1}$.

In recent years, fractional calculus has attracted a lot of mathematicians' interests.  The CR fractional sub-Laplacian $P_\gamma^{\theta_1}(\gamma \in (0, 1))$ is defined by Rupert L. Frank, Mari a del~Mar Gonzalez, Dario D. Monticelli, and Jinggang Tan in \cite{FrankGonzalezMonticelliTan}. In the paper \cite{FrankGonzalezMonticelliTan}, it was shown that one can  treat the CR fractional sub-Laplacian as a boundary operator. In \cite{BFM2013}, the CR fractional sub-Laplacian is viewed as intertwining operator.

On the CR sphere, the general intertwining operator $P_\gamma^{\theta_1}$ is defined by the
following property:
\begin{align*}
|J_\tau|^{\frac{n+1+\gamma}{2n+2}}(P_\gamma^{\theta_1}\circ F)\circ \tau=P_\gamma^{\theta_1}\left( |J_\tau|^{\frac{n+1-\gamma}{2n+2}}( F\circ \tau)\right), \quad \forall \tau\in \textmd{Aut}(\mathbb{S}^{2n+1})
\end{align*}
for each $F\in C^{\infty}(\mathbb{S}^{2n+1})$.

Given another conformal representative $\tilde{\theta}=u^{\frac{2}{n+1-\gamma}}\theta_1$
which identifies the CR structure of $\mathbb{S}^{2n+1}$, the corresponding operator is defined by
$$P_\gamma^{\tilde{\theta}}(\phi)=u^{-\frac{n+1+\gamma}{n+1-\gamma}}P_\gamma^{\theta_1}(u\phi).$$

For $(\mathbb{S}^{2n+1}, \theta)$, we define its fractional Q-curvature as
$$Q_\gamma^{\theta}=P_\gamma^{\theta}(1).$$
Let $\phi=1$,
$$P_\gamma^{\theta_1}(u)=u^{\frac{n+1+\gamma}{n+1-\gamma}}P_\gamma^{\tilde{\theta}}(1).$$

The
fractional Q-curvature problem is that for a prescribed  function $K$, whether there is a positive function $u$ such that $Q_\gamma^{\tilde\theta}=K$ with $\tilde{\theta}=u^{\frac{2}{n+1-\gamma}}\theta_1$. This problem is equivalent to the
 existence of  the following  fractional nonlinear PDE:
\begin{equation}\label{main}
P_\gamma^{\theta_1}u=K u^{\frac{n+1+\gamma}{n+1-\gamma}}, \;u>0 \ \textmd{on} \ \ \mathbb{S}^{2n+1}.
\end{equation}

The scalar curvature problem for the Riemannian manifolds  has been extensively studied, in dimension 2,3 (see \cite{bahri1991},\cite{ChangYang1987},\cite{hz90}) as well as in high dimension (\cite{MBA1996},\cite{lyy95}).
Fractional scalar curvature problem for the Riemannian manifolds has been studied by \cite{cz14},\cite{cz15}. There are also many works on scalar curvature problem on CR mainfolds , see \cite{ChtiouiEl2007},\cite{GamaraNajoua2001},\cite{3CR02},\cite{MalchiodiUguzzoniJMPA2002},\cite{sa10}. In this paper, we study
fractional Q-curvature problem on $\mathbb{S}^{3}$.

 Let
$\dot{\mathcal{S}}^{\gamma}(\mathbb{S}^{2n+1})$ be the closure of $\mathcal{C}_0^\infty(\mathbb{S}^{2n+1})$ with respect to the quadratic form
$$
\int_{\mathbb{S}^{2n+1}}fP_\gamma^{\theta_1} f\theta_1\wedge(d\theta_1)^n,
$$
$\Sigma=\{u\in\dot{\mathcal{S}}^{\gamma}(\mathbb{S}^{2n+1})|  \|u\|=\sqrt{\int_{\mathbb{S}^{2n+1}}uP_\gamma u\theta_1\wedge(d\theta_1)^n}=1\}$ and
$\Sigma^+=\{u\in\Sigma  | u\geqslant 0\}$.

For $u\in\dot{\mathcal{S}}^{\gamma}(\mathbb{S}^{2n+1})$, we define
$$J(u)=\frac{\|u\|^2}{(\int_{\mathbb{S}^{2n+1}}K u^{\frac{Q}{n+1-\gamma}}\theta_1\wedge(d\theta_1)^n)^{\frac{n+1-\gamma}{n+1}}}.$$
If $u$ is a critical point of the function $J$ in $\Sigma^+$, then $v=(J(u))^{\frac{n+1}{2\gamma}}u$ is a solution of (\ref{main}).
However, the functional $J$ does not satisfy the Palais-Smale condition, that is to say there exist
critical points at infinity, which are the limits of noncompact orbits for
the gradient flow of $-J$. Thinking of these sequences as critical points, a natural idea is to expand the functional
$J$ near the sets of such critical points.

In this paper we care the case $n=1$.
We state now the main result.
If $K: \mathbb{S}^{3}\rightarrow \mathbb{R}$ is a $C^2$ positive function, we assume $K$ satisfying condition:
\begin{align}\label{nd}
\textmd{each critical ponit} \ \eta_i \ \textmd{is a non degenerate critical ponit of}\ K \ \textmd{and}\ \Delta_{\theta_1} K(\eta_i)\neq 0.
\end{align}
Denote
$$I^+=\{\eta_i\in \mathbb{S}^{3}: \nabla_{\theta_1}  K(\eta_i)=0, \ -\Delta_{\theta_1} K(\eta_i)> 0\}.$$
Assume that $\sharp I^+=m$, $\sharp I^+$ is the cardinality of $I^+$. For simplicity, we assume $I^+=\{\eta_1, \eta_2,\cdots, \eta_{m}\}$.
For any $l$-element subset $\{\eta_{i_1},\cdots, \eta_{i_l}\}$ of $I^+$, $1\le l\le m$, we define $\mu_{\tau_l}=\sum_{j=1}^l\textmd{ind}(K,\eta_{i_j})$ with $\tau_l=(i_1,\cdots,i_l)$.
\begin{theorem}\label{the}
Let $\frac{2}{3}\leqslant\gamma<1$, assume $K$ satisfies (\ref{nd}). Then the problem (\ref{main}) has a solution provided
$$\sum_{l=1}^{m}\sum_{\tau_l}(-1)^{\mu_{\tau_l}}\neq -1.$$
\end{theorem}

We will prove the theorem by contradiction in section 5. Therefore we assume that equation
(\ref{main}) has no solutions. Our proof is based on a  technical Morse Lemma at
infinity; it relies the construction of a suitable pseudogradient for
$J$. The (PS) condition is satisfied along the decreasing flow lines of
this pseudogradient, as long as these flow lines stay out of the neighborhood of a finite
number of critical points of $K$. Finally we compute the contribution of some critical points  at infinity to the
changes of topology for   the level set of the functional, from this we can achieve  a contradiction.

This paper is organized as follows. In the next section, we introduce preliminary result and the general variational
framework . In section 3, we give some expansions of the functional and its gradient
near the sets of its critical points at infinity. In section 4, we establish the Morse lemma
at infinity, which allows us to refine the expansion of the function. In section 5, we give
a proof of Theorem \ref{the}. In Appendices
A-C, we  show some useful estimates which will be used in our proof of Theorem \ref{the}.

\section{Preliminary results}
The Heisenberg group $\mathbb{H}^1$ is a Lie group whose underlying manifold is $\mathbb{R}\times\mathbb{R}\times\mathbb{R}$ with elements $u = (x, y, t)$ and whose group law is
$$
u\circ u' = (x, y, t)\circ(x', y', t') = \left(x + x', y + y',  t + t' + 2(x'y-xy')\right).
$$

Alternatively, we can use complex coordinates $z=x+iy$ to denote elements of $\mathbb{R}\times\mathbb{R}\simeq \mathbb{C}$, so that the group law can be written as
$$(z',t')\circ (z,t)=(z'+z, t'+t+2\textmd{Im}<z',z>_{\mathbb{C}}),$$
for  $(z',t'), (z,t)\in \mathbb{H}^1$, and $<\cdot,\cdot>_{\mathbb{C}}$ is the standard Hermitian inner product in $\mathbb{C}$.

The  CR sphere $\mathbb{S}^{2n+1}=\{\zeta=(\zeta_1,\cdots,\zeta_{n+1})\in \mathbb{C}^{n+1}, | \sum_{j=1}^{n+1}|\zeta_j|^2=1\}.$ The standard Euclidean volume
element of $\mathbb{S}^{2n+1}$ is denoted by $\textmd{d}\zeta$.

We introduce Cayley transform $\mathcal{C}$ between the Heisenberg group and the  CR sphere.
\begin{align*}
\mathcal{C}: \mathbb{H}^{n} &\rightarrow \mathbb{S}^{2n+1}\setminus(0,0,\cdots,-1),\\
(z,t) &\mapsto (\frac{2z}{1+|z|^2+it}, \frac{1+|z|^2-it}{1+|z|^2+it}).
\end{align*}
The inverse is given by
\begin{align*}
\mathcal{C}^{-1}: \mathbb{S}^{2n+1}\setminus(0,0,\cdots,-1) &\rightarrow \mathbb{H}^{n}, \\ \zeta=(\zeta_1,\cdots,\zeta_{n+1}) &\mapsto (\frac{\zeta_1}{1+\zeta_{n+1}},\cdots,\frac{\zeta_n}{1+\zeta_{n+1}},\textmd{Im}\frac{1-\zeta_{n+1}}{1+\zeta_{n+1}}).
\end{align*}
The Jacobian determinant of $\mathcal{C}$ is
$$|Jac(z,t)|=\frac{2^{2n+1}}{((1+|z|^2)^2+t^2)^{n+1}}.$$

For any $\lambda > 0$ the dilation $\lambda:\mathbb{H}^1\to \mathbb{H}^1$ is defined by $\lambda u = \lambda (x, y, t) = (\lambda x,\lambda y, \lambda^2 t)$ and we denote the homogeneous norm on $\mathbb{H}^1$ by $|u| = |(x, y, t)| = ((x^2+y^2)^2 + t^2)^{1/4}$.  The CR structure on $\mathbb{H}^1$ is given by the left invariant vector field:
$$X=\frac{\partial}{\partial x}+2y\frac{\partial}{\partial t}, \;\;\; Y=\frac{\partial}{\partial y}-2x\frac{\partial}{\partial t}, \;\;\;
T=\frac{\partial}{\partial t}.$$
The standard contact form $\theta_0=\textmd{d}t+2(x\textmd{d}y-y\textmd{d}x)$.  Haar measure on $\mathbb{H}^1$ is the Lebesgue measure $du = dxdydt$. Denote $\nabla_{\theta_0}=(X,Y,T)$.
In $\mathbb{H}^1$, Taylor polynomials can be written in a special symmetric form. The expansions are similar to Taylor expansions in $\mathbb{R}^n$ but are adjusted to compensate for the different Heisenberg structure.
The following formula from \cite{BieskeThomas2002} gives the Taylor expansions based at the origin.
Let $f: \mathbb{H}^1\rightarrow \mathbb{R}$ be a $C^2$ functions. Let the origin be denoted by $0$, $p=(x,y,t)$ be an arbitrary point around $0$. Then,
\begin{align*}
f(p)&=f(0)+x(Xf)(0)+y(Yf)(0)+t(Tf)(0)+\frac{x^2}{2}(X^2f)(0)+\frac{y^2}{2}(Y^2f)(0)\\&\quad+\frac{xy}{2}(XYf)(0)+\frac{xy}{2}(YXf)(0)+o(|p|^2).
\end{align*}
The sub-Laplacian on $\mathbb{H}^1$ is the second order differential operator
$$\Delta_{\theta_0}=\frac{1}{4}(X^2+Y^2).$$

On the CR sphere, the standard contact form $\theta_1=i\sum_{j=1}^{n+1}(\zeta_j\textmd{d}\bar{\zeta}_j-\bar{\zeta}_j\textmd{d}\zeta_j)$, the sub-gradient is $\nabla_{\theta_1}$ and the sub-Laplacian is defined as
$$\Delta_{\theta_1}=\frac{1}{2}\sum_{j=1}^{n+1}(T_j\bar{T}_j+\bar{T}_jT_j),$$
where
$$T_j=\frac{\partial}{\partial \zeta_j}-\bar{\zeta}_jR,  \quad  R=  \sum_{k=1}^{n+1} \zeta_k\frac{\partial}{\partial \zeta_k}.$$

The conformal sub-Laplacian on the sphere is defined as
$$L_{\theta_1}=\Delta_{\theta_1}+\frac{1}{4}.$$
The peculiarity of $L_{\theta_1}$ is its direct relation with $\Delta_{\theta_0}$ via the Cayley transform:
\begin{align}\label{relat1}
\Delta_{\theta_0}\left((2|Jac|)^{\frac{1}{4}}(F\circ \mathcal{C})\right)=(2|Jac|)^{\frac{3}{4}}(L_{\theta_1}F)\circ\mathcal{C},
\end{align}
where $F:\mathbb{S}^{2n+1}\rightarrow \mathbb{C}$ is a smooth  function.

The differences between the standard volume elements for $\mathbb{S}^{2n+1}$ and $\mathbb{H}^{n}$ and the volume forms associated
with the standard contact forms $\theta_1$, and $\theta_0$ of these two spaces state as:
\begin{align}\label{relat2}
\int_{\mathbb{S}^{2n+1}}F \theta_1\wedge (d\theta_1)^n&=2^{2n+1}n!\int_{\mathbb{S}^{2n+1}}F\textmd{ d}\zeta\nonumber\\&=\int_{\mathbb{H}^{n}}2|Jac|F\circ \mathcal{C} \theta_0\wedge(d\theta_0)^n\nonumber\\&=2^{2n}n!\int_{\mathbb{H}^{n}}2|Jac|F\circ \mathcal{C}\textmd{d}u.
\end{align}
We refer \cite{BFM2013} for details.

We consider the CR fractional operators of order $2\gamma$. For $\gamma \in (0, 1)$, the symbol of the operator on $\mathbb{H}^n$ is
\begin{equation*}
P_\gamma^{\theta_0} =(2|T|)^\gamma\displaystyle\frac{\Gamma\left(\frac{1 + \gamma}{2} + \frac{\Delta_{\theta_0}}{2|T|}\right)}{\Gamma\left(\frac{1 - \gamma}{2} + \frac{\Delta_{\theta_0}}{2|T|}\right)}.
\end{equation*}
In particular, for $\gamma=1$, $P_1^{\theta_0}=\Delta_{b}$ is the CR Yamabe operator on  $\mathbb{H}^n$.

$P_\gamma^{\theta_0}$ also satisfy (setting $Q=2n+2$ and $d=2\gamma$)
\begin{align*}
|J_h|^{\frac{Q+d}{2Q}}(P_\gamma^{\theta_0}\circ f)\circ h=P_\gamma^{\theta_0}( |J_h|^{\frac{Q-d}{2Q}}( f\circ h)), \quad \forall h\in \textmd{Aut}(\mathbb{H}^{n}),
\end{align*}
for each $f\in C^{\infty}(\mathbb{H}^{n}).$

On the CR sphere, the general intertwining operator $P_\gamma^{\theta_1}$ is defined by the
following property:
\begin{align*}
|J_\tau|^{\frac{Q+d}{2Q}}(P_\gamma^{\theta_1}\circ F)\circ \tau=P_\gamma^{\theta_1}\left( |J_\tau|^{\frac{Q-d}{2Q}}( F\circ \tau)\right), \quad \forall \tau\in \textmd{Aut}(\mathbb{S}^{2n+1}),
\end{align*}
for each $F\in C^{\infty}(\mathbb{S}^{2n+1}).$

Thus we have
\begin{align}\label{relat3}
P_\gamma^{\theta_0}\left((2|Jac|)^{\frac{Q-d}{2Q}}(F\circ \mathcal{C})\right)=(2|Jac|)^{\frac{Q+d}{2Q}}(P_\gamma^{\theta_1}F)\circ\mathcal{C}.
\end{align}
For more details, we refer to the well written paper \cite{FrankGonzalezMonticelliTan} and \cite{BFM2013}.

It was proved by
\cite{FrankLiebAnnals2013} that on $\mathbb{H}^n$:
for $q = \frac{Q}{n+1- \gamma}$ and any function $f\in \dot{\mathcal{S}}^{\gamma}(\mathbb{H}^n)$, it holds that
\begin{equation}\label{e:FLHLS1}
\|f\|_{L^q(\mathbb{H}^n)}^2\leq \frac{2^{n - \gamma - 1 - \frac{2(n - 1)\gamma}{Q}}(\pi^{n + 1})^{2\gamma/Q - 1}}{(n!)^{2\gamma/Q - 1}}\frac{\Gamma((Q - 2\gamma)/2)\Gamma^2((Q - 2\gamma)/4)}{\Gamma^2((Q - \gamma)/2)\Gamma(\gamma)}\int_{\mathbb{H}^n}fP_\gamma fdu
\end{equation}
and all optimizers are translates, dilates or constant multiples of the function
$$\delta(u) = \delta(z, t) = \left(\frac{1}{(1 + |z|^2)^2 + t^2}\right)^{\frac{Q - 2\gamma}{4}}.$$

We know that, for $\lambda > 0$, $a\in\mathbb{H}^1$ and some suitable choice of $c_0 = C(\gamma) > 0$, the function
\begin{equation}\label{e:extremals}
\delta_{a, \lambda}(u) = c_0\lambda^{2-\gamma} \delta(\lambda(a^{-1}u))
\end{equation}
satisfies the Euler-Lagrangian equation for
\begin{equation}\label{e:CRYamabe}
P_\gamma^{\theta_0} u = u^{\frac{2+\gamma}{2- \gamma}},\quad u > 0 \text{ in } \mathbb{H}^1.
\end{equation}
(\ref{e:CRYamabe}) is the fractional CR Yamabe equation introduced in \cite{FrankGonzalezMonticelliTan}.
(\ref{e:extremals}) indicates that (\ref{e:CRYamabe}) is invariant under the scaling and translations.

We introduce the function for each $(\zeta_0,\lambda)\in \mathbb{S}^{2n+1}\times (0,+\infty)$,
\begin{align}
w_{\zeta_0,\lambda}(\zeta)=|1+\zeta_{n+1}|^{-(2-\gamma)}\delta_{\mathcal{C}^{-1}(\zeta_0), \lambda}\circ\mathcal{C}^{-1}(\zeta).
\end{align}

Using (\ref{relat2}) and (\ref{relat3}), we have
\begin{align*}
\int_{\mathbb{S}^{3}}w_{\zeta_0,\lambda}^{\frac{2-\gamma}{4}}\theta_1\wedge d\theta_1=\int_{\mathbb{H}^1}\delta_{a_0, \lambda}^{\frac{2-\gamma}{4}}
\theta_0\wedge d\theta_0,
\end{align*}
and
\begin{align*}
\int_{\mathbb{S}^{3}}w_{\zeta_0,\lambda}P_\gamma^{\theta_1}w_{\zeta_0,\lambda}\theta_1\wedge d\theta_1=\int_{\mathbb{H}^1}\delta_{a_0, \lambda}P_\gamma^{\theta_0}\delta_{a_0, \lambda}
\theta_0\wedge d\theta_0,
\end{align*}
where $a_0=\mathcal{C}^{-1}(\zeta_0),\xi=\mathcal{C}^{-1}(\zeta)$. We also have $P_\gamma^{\theta_1}w_{\zeta_0,\lambda}=w_{\zeta_0,\lambda}^{\frac{2+\gamma}{2-\gamma}}$.

For any $\varepsilon >0$ and $p \geqslant 1$, we set
$$(\alpha,g,\lambda)=(\alpha_1,\cdots,\alpha_p,g_1,\cdots,g_p,\lambda_1,\cdots,\lambda_p)\in (0,+\infty)^p\times(\mathbb{S}^{3})^p\times(0,+\infty)^p,$$
$$\varepsilon_{ij}=\left(\frac{\lambda_i}{\lambda_j}+\frac{\lambda_j}{\lambda_i}+\lambda_i\lambda_jd(g_i,g_j)^2\right)^{-(2-\gamma)}.$$
Let $V(p,\varepsilon)$ be the subset of $\Sigma^{+}$ of the following functions: $u\in\Sigma^+,\exists (\alpha,g,\lambda),$ such that
$$\|u-\sum_{i=1}^p\alpha_i w_{g_i,\lambda_i}\|<\varepsilon$$
and $|J(u)^{\frac{2}{2-\gamma}}\alpha_i^{\frac{2\gamma}{2-\gamma}}K(g_i)-1|<\varepsilon$, $\varepsilon_{ij}<\varepsilon$,$\lambda_i>\varepsilon^{-1}$.
The set $V(p,\varepsilon)$ has a simple interpretation: It is a neighborhood of the
critical points at infinity of the functional $J$ on $\Sigma^+$.

\begin{definition}\cite{babook}
We will say that the Palais-Smale condition holds on flow-lines in the $V(p,\varepsilon)$ if, taking an initial data $u_0$ in $V(p,\varepsilon)$, with $\varepsilon_0$ small enough (but fixed), the solution $u(s,u_0)$ of the differential equation $\frac{\partial u}{\partial s}=-\partial J(u)$ with initial data $u_0$ remains outside a $V(p,\varepsilon_1), \varepsilon_1>0$, which depends only on $u_0$.
\end{definition}

The failure of (PS)condition is characterized as follows.
\begin{lemma}\label{PSfail}
Assume that (\ref{main}) has no solutions. Let $\{u_k\}\subseteq\Sigma^{+}$ be a sequence such that $J'(u_k)\rightarrow 0$ and $J(u_k)$ is bounded. Then
there exists an integer $p \geqslant 1$, a positive sequence $\varepsilon_k \rightarrow0$ and an extracted subsequence of $\{u_k\}$, such that $u_k\in
V(p,\varepsilon_k)$.
\end{lemma}

In order to prove Lemma \ref{PSfail}, we introduce
$$I(u)=\frac{1}{2}\|u\|^2-\frac{2-\gamma}{4}\int_{\mathbb{S}^{3}}K u^{\frac{4}{2-\gamma}}\theta_1\wedge d\theta_1.$$
Note that $J'(u_k)\rightarrow 0$ if and only if $I'(J(u_k)^{\frac{1}{\gamma}}u_k)\rightarrow 0$. Then we can follow the first part of \cite{coron} by using the functional $I$. The proof is by now classical, we can also see \cite{gamara01}.

We introduce the minimization problem for $\varepsilon$ small enough
\begin{equation}\label{miniza}
 \displaystyle \textmd{min}_{u\in V(p,\varepsilon)}\{\|u-\sum_{i=1}^p\alpha_iw_{g_i,\lambda_i}\|, \ \alpha_i>0, g_i\in\mathbb{S}^{3}, \lambda_i >0\}.
\end{equation}
\begin{lemma}\label{Vv}
For any $p \geqslant 1$, there exists $\varepsilon_p>0$ such that for any $0<\varepsilon<\varepsilon_p$ and $u\in V(p,\varepsilon)$, the minimization problem
(\ref{miniza}) has a unique solution $(\bar{\alpha},\bar{g},\bar{\lambda})$. Denoting $v=u-\sum_{i=1}^p\alpha_iw_{g_i,\lambda_i}$,$v$
satisfies
\begin{equation}\label{V0}
\begin{cases}
\langle v,w_{g_i,\lambda_i}\rangle=0\\ \langle v, \frac{\partial w_{g_i,\lambda_i}}{\partial\lambda_i}\rangle=0\\
\langle v, \frac{\partial w_{g_i,\lambda_i}}{\partial g_i}\rangle=0 \ i=1,\cdots,p.
\end{cases}
\end{equation}
Here $\langle \cdot,\cdot\rangle$ denote the inner product in $\dot{\mathcal{S}}^{\gamma}(\mathbb{S}^{3})$ defined by
$$\langle u,v\rangle=\int_{\mathbb{S}^{3}}uP_\gamma^{\theta
_1} v\theta_1\wedge d\theta_1.$$
\end{lemma}
The proof of Lemma \ref{Vv} follows from Appendix A in \cite{BahriCoron1988} with some modulations.

\section{The expansion of the function $J$}
\begin{lemma}\label{expansion}
If $\varepsilon >0$ small enough and $u=\sum_{i=1}^p\alpha_i w_{g_i,\lambda_i}+v\in V(p,\varepsilon)$, $v$ satisfies (\ref{V0}), we have
\begin{align*}
J(u)&=\frac{\sum_{i=1}^p\alpha_i^2S}{\left(\sum_{i=1}^p\alpha_i^{\frac{4}{2-\gamma}}K(g_i)\right)^{\frac{2-\gamma}{2}}}\left[1-\frac{2-\gamma}{2}\frac{c_2}{S^{\frac{2}
{\gamma}}}\sum_{i=1}^p\frac{\alpha_i^{\frac{4}{2-\gamma}}\Delta_{\theta_1} K(g_i)}{\sum_{k=1}^p\alpha_k^{\frac{4}{2-\gamma}}K(g_k)\lambda_i^2}\right.\\&\quad
\left.+(f,v)+Q(v,v)+O'(\sum_{i\neq j} \varepsilon_{ij})+o(\sum_{i=1}^p\frac{1}{\lambda_i^2})+o(\|v\|^2)\right],
\end{align*}
with
\begin{align*}
(f,v)=\frac{-2}{\sum_{k=1}^p\alpha_k^{\frac{4}{2-\gamma}}K(g_k)S^{\frac{2}{\gamma}}}\int_{\mathbb{S}^3}K(\zeta)\left(\sum_{i=1}^p\alpha_iw_{g_i,\lambda_i}\right)^
{\frac{2+\gamma}{2-\gamma}}v\theta_1\wedge d\theta_1,
\end{align*}
\begin{align*}
Q(v,v)&=\frac{1}{\sum_{k=1}^p\alpha_k^2S^{\frac{2}{\gamma}}}\|v\|^2\\&\quad-\frac{2+\gamma}{(2-\gamma)S^{\frac{2}{\gamma}}\sum_{k=1}^p\alpha_k^{\frac{4}
{2-\gamma}}K(g_k)}\int_{\mathbb{S}^3}K(\zeta)\sum_{i=1}^p(\alpha_iw_{g_i,\lambda_i})^{\frac{2\gamma}{2-\gamma}}v^2.
\end{align*}
Furthermore $\|f\|$ is bounded by
\begin{align*}
\|f\|= O\left(\sum_{i=1}^p\left(\frac{|\nabla_{\theta_1}K(g_i)|}{\lambda_i}+\frac{1}{\lambda_i^2}\right)+\sum_{i\neq j}\varepsilon_{ij}(\ln\varepsilon_{ij}^{-1})^{{\frac{2-\gamma}{2}}}\right).
\end{align*}
\end{lemma}

The proof of this lemma is provided in Appendix A.

Now, we state the following two lemmas whose proof follow the arguments used to prove similar statements in \cite{bahri1991} (also in \cite{babook}); see the Appendix of \cite{3CR02} where some necessary modifications
are made.
\begin{lemma}
$Q(v,v)$ is a quadratic form positively definite in
$$H_\varepsilon(\lambda, a)=\{v\in\dot{\mathcal{S}}^{\gamma}(\mathbb{S}^{3})| \ \ v\,\textmd{satisfies }(\ref{V0}),\ \|v\|\leqslant\varepsilon \}.$$
\end{lemma}
One can follow the idea of the proof of Lemma A.2 in \cite{bahri1991} to get a proof of this lemma. We omit the details.
\begin{lemma}\label{vbar}
For any $u_0=\sum_{i=1}^p\alpha_iw_{g_i,\lambda_i}\in V(p,\varepsilon)$, there exists a unique $\overline{v}=\overline{v}(\alpha, g, \lambda)$
which minimizes $J(u_0+v)$ with respect to $v\in H_\varepsilon(\lambda, a) $ and we have estimate $\|\overline{v}\|= O(\|f\|)$.
\end{lemma}

For a proof of Lemma \ref{vbar}, one may follow the idea and similar estimates in the proof of Proposition 5.4 in \cite{babook}(P191). We omit the details.

Since $\overline{v}$ is a minimizer, we have
$$(f,\overline{v})+2Q(\overline{v},\overline{v})+o(\|\overline{v}\|^2)=0.$$
It yields
$$(f,v)+Q(v,v)+o(\|v\|^2)=Q(v-\overline{v},v-\overline{v})-Q(\overline{v},\overline{v})+o(\|\overline{v}\|^2).$$
From Lemma 3.1 and Lemma 3.3, we state the following lemma which improve the asymptotic behavior of the function $J$.
\begin{lemma}\label{expansion1}
For any $p\geqslant 1$, there exists $\varepsilon_p >0$ such that for $u=\sum_{i=1}^p\alpha_i w_{g_i,\lambda_i}+v$, $v\in H_\varepsilon(\lambda, a)$, we have
\begin{align*}
J(u)&=\frac{\sum_{i=1}^p\alpha_i^2S}{\left(\sum_{i=1}^p\alpha_i^{\frac{4}{2-\gamma}}K(g_i)\right)^{\frac{2-\gamma}{2}}}\left[1-\frac{2-\gamma}{2}\frac{c_2}{S^{\frac{2}
{\gamma}}}\sum_{i=1}^p\frac{\alpha_i^{\frac{4}{2-\gamma}}\Delta_{\theta_1} K(g_i)}{\sum_{k=1}^p\alpha_k^{\frac{4}{2-\gamma}}K(g_k)\lambda_i^2}\right.\\&\quad\left.
+Q(v-\overline{v},v-\overline{v})-Q(\overline{v},\overline{v})+o(\|\overline{v}\|^2)+O(\sum_{i\neq j} \varepsilon_{ij})+o(\sum_{i=1}^p\frac{1}{\lambda_i^2})\right].
\end{align*}
\end{lemma}

\begin{lemma}
Assume $K$ be a $C^2$ positive function satisfying condition (\ref{nd}). For any $u=\sum_{i=1}^p\alpha_iw_{g_i,\lambda_i}\in V(p,\varepsilon)$, the following estimate holds
\begin{align}\label{es01}
J'(u)\left(\lambda_j\frac{\partial w_{g_j,\lambda_j}}{\partial\lambda_j}\right)&= \nonumber2\lambda(u)\left[-\sum_{i\neq j}O'(\lambda_j\frac{\partial\varepsilon_{ij}}{\partial\lambda_j})+
\frac{2-\gamma}{4}c_2\alpha_j\frac{\Delta_{\theta_1} K(g_j)}{ K(g_j)\lambda_j^2}\right.\\&\quad\left.o(\frac{1}{\lambda_j^2})+
o(\sum_{i\neq j} \varepsilon_{ij})\right],
\end{align}
\begin{align}\label{es02}
J'(u)\left(\frac{1}{\lambda_j}\frac{\partial w_{g_j,\lambda_j}}{\partial g_j}\right)=-2\lambda(u)\alpha_j
c_2\frac{\nabla_{\theta_1}K(g_j)}{ K(g_j)\lambda_j}+
O\left(\sum_{i\neq j} \varepsilon_{ij}+\frac{1}{\lambda_j^2}\right).
\end{align}
\end{lemma}
We will give the proof in Appendix B and Appendix C.

\section{Morse lemma at infinity}

This section is devoted to characterize the critical points at infinity associated to problem (\ref{main}). The characterization is obtained through the construction of a suitable pseudogradient at infinity in the set $V(p,\varepsilon)$, depending on a delicate expansion of the gradient of J near infinity.

\begin{theorem}\label{covering}
There is a covering $\{O_l\}$ and a subset of $\{(\alpha_l,g_l,\lambda_l)\}$ of the base space for the
bundle $V(p,\varepsilon)$ and a diffeomorphism $\xi_l : V(p,\varepsilon)\rightarrow V(p,\varepsilon')$ for some $\varepsilon'> 0$ with
$$\xi_l\left(\sum_{i=1}^p\alpha_i w_{g_i,\lambda_i}+\overline{v}(\alpha,g,\lambda)\right)=\sum_{i=1}^p\alpha_iw_{\tilde{g_i},\tilde{\lambda_i}}$$ such that
\begin{align}\label{zz}J\left(\sum_{i=1}^p\alpha_iw_{g_i,\lambda_i}+v\right)=J\left(\sum_{i=1}^p\alpha_iw_{\tilde{g_i},\tilde{\lambda_i}}\right)+\frac{1}{2}J''
\left(\sum_{i=1}^p\alpha_iw_{
{g_i},{\lambda_i}}+\bar v(\alpha,g,\lambda)\right)V_l\cdot V_l,\end{align}
where $(\alpha,g,\lambda)\in O_l$, $(\alpha,\tilde{g},\tilde{\lambda})$ not depending on $O_l$, $V_l$ orthogonal to $w_{\tilde{g_i},\tilde{\lambda_i}},\frac{\partial w_{\tilde{g_i},\tilde{\lambda_i}}}{\partial \tilde{\lambda_i}},\frac{\partial w_{\tilde{g_i},\tilde{\lambda_i}}}{\partial \tilde{g_i}}$.
\end{theorem}
The proof of this theorem need some technical result. First we give the Morse lemma at infinity by isolating the contribution of $v-\overline{v}$.
\begin{lemma}\label{4.1}
For any $\sum_{i=1}^p\alpha_i w_{\bar{g_i},\bar{\lambda_i}}\in V(p,\varepsilon)$, let
$$(\bar{\alpha},\bar{g},\bar{\lambda})=(\bar{\alpha_1},\cdots,\bar{\alpha_p},\bar{g_1},\cdots,\bar{g_p},\bar{\lambda_1},
\cdots,\bar{\lambda_p}),$$
there is a neighborhood $U$ of $(\bar{\alpha},\bar{g},\bar{\lambda})$ such that
$$J\left(\sum_{i=1}^p\alpha_iw_{g_i,\lambda_i}+v\right)=J\left(\sum_{i=1}^p\alpha_iw_{g_i,\lambda_i}+\bar{v}(\alpha,g,\lambda)\right)+\frac{1}{2}J''
\left(\sum_{i=1}^p
\bar{\alpha_i}w_{\bar{g_i},\bar{\lambda_i}}+\bar{v}(\bar{\alpha},\bar{g},\bar{\lambda})\right)V\cdot V$$
for any $u=\sum_{i=1}^p\alpha_iw_{g_i,\lambda_i}+v\in V(p,\varepsilon)$ with $(\alpha,g,\lambda)\in U$, where $V=V(\alpha,g,\lambda,v)$ is a $C^1$
diffeomorphism that has range orthogonal to
$$\bigcup_{i=1}^p\left\{\delta_{{g_i'},{\lambda_i'}},\frac{\partial\delta_{{g_i'},{\lambda_i'}}}{\partial {\lambda_i'}},\frac{\partial\delta_{{g_i'},{\lambda_i'}}}{\partial {g_i'}}\right\}$$
for any $(\alpha',g',\lambda')\in U$ and $\|V\|=O(\|v-\bar{v}(\bar{\alpha},\bar{g},\bar{\lambda})\|)$.
\end{lemma}
The proof of Lemma \ref{4.1} is similar to the proof of Lemma 3.2  in \cite{MBA1996} for the Riemannian manifold.
\begin{lemma}\label{Wcon}
Let $K\in C^2(\mathbb{S}^3)$ be a positive function satisfying condition (\ref{nd}). For any $u=\sum_{i=1}^p\alpha_iw_{g_i,\lambda_i}\in V(p,\varepsilon), \varepsilon$ small enough, Then there exists a vector filed $W'$ so that the following holds: there is a constant
$C>0$  such that
\begin{align}\label{w1}
-J'(u)(W')\geqslant C_1\left(\sum_{i=1}^p(\frac{|\nabla_{\theta_1}K(g_i)|}{\lambda_i}+\frac{1}{\lambda_i^2})+\sum_{i\neq j} \varepsilon_{ij}\right),\end{align}
\begin{align}\label{w2}
-J'(u+\overline{v})\left(W'+\frac{\partial\overline{v}}{\partial(\alpha,g,\lambda)}(W')\right)\geqslant C_1\left(\sum_{i=1}^p(\frac{|\nabla_{\theta_1} K(g_i)|}{\lambda_i}+\frac{1}{\lambda_i^2})+\sum_{i\neq j} \varepsilon_{ij}\right),
\end{align}
where $C_1$ is a positive constant, $\|W'\|$ is bounded.
\end{lemma}

We remind that $\frac{\partial\overline{v}}{\partial(\alpha,g,\lambda)}(W)=\frac{\partial\overline{v}}{\partial W}$ means the variation of $\bar v$ along the direction  $W$, where $W=W_{\alpha}\delta \alpha+W_{g}\cdot \delta g+W_{\lambda}\delta\lambda$ is an increment in the $(\alpha_i,g_i,\lambda_i)$ space.

 We define a set $I_i$ for all $i\in \{1,\cdots,p\}$. We   divide it in three  cases:

\begin{enumerate}
\item  For all $i\in\{1,\cdots,p\}$,  there exists  a suitable constant $C>0$ such that
\begin{align}\label{case1}
\sum_{j \neq i} \varepsilon_{ij}\leqslant \frac{C}{\lambda_i^2}.
\end{align}

In this case, we define $I_i=\emptyset$ the empty set.

\item   $\lambda_i$ is the largest concentration with
\begin{align}\label{case2}
\sum_{j\neq i} \varepsilon_{ij}> \frac{C}{\lambda_i^2},
\end{align} then in this case we define $I_i=\{i\}$.

\item  $\lambda_i$ satisfies \eqref{case2}, but it is not the largest concentration, {\it i.e.}, we have another $j$ such that $\lambda_j>\lambda_i$ and $\lambda_j$ satisfies \eqref{case2}. In this case we define
\begin{align}\label{case3p}I_i=\left\{k \ | \ \lambda_k\geqslant\lambda_i, \ \sum_{j\neq k} \varepsilon_{kj}> \frac{C}{\lambda_k^2} \right\}.\end{align}

\end{enumerate}

Suppose $I_i=\{k_1, k_2,\cdots, k_m\}$ such that $k_1=i$ and $\lambda_{k_1}\le \lambda_{k_2}\le \cdots \le\lambda_{k_m}$. We define $\mu_{k_s}=2^{s-1}$, $s=1,\cdots,m$.

\noindent{\it Claim.} There holds
\begin{align}\label{esaa}
\left|J'(u)\left(\frac{1}{\lambda_i}\frac{\partial w_{g_i,\lambda_i}}{\partial g_i}\right)\right|+\frac{c'}{c}\sum_{k\in I_i}\mu_{k}J'(u)\left(\lambda_{k}\frac{\partial w_{g_{k},\lambda_{k}}}{\partial\lambda_{k}}\right)\geqslant c\frac{|\nabla_{\theta_1} K(g_i)|}{\lambda_i}-\frac{1}{c}\frac{1}{\lambda_i^2},
\end{align} where $c$ is the constant in \eqref{es02} and $c'>0$ is a suitable constant.

\noindent{\it Proof of the claim}. In the case (1),  by \eqref{es02}, we have
\begin{align}\label{estimate1}
\left|J'(u)\left(\frac{1}{\lambda_i}\frac{\partial w_{g_i,\lambda_i}}{\partial g_i}\right)\right|\geqslant
c\frac{|\nabla_{\theta_1}K(g_i)|}{\lambda_i}-\frac{1}{c}\left(\sum_{j\neq i} \varepsilon_{ij}+\frac{1}{\lambda_i^2}\right)\geqslant c\frac{|\nabla_{\theta_1} K(g_i)|}{\lambda_i}-\frac{1}{c}\frac{1}{\lambda_i^2}.
\end{align}

In the case (2),
then for $\lambda_j>\lambda_i$, there holds
\begin{align}\label{es21}
\sum_{k\neq j} \varepsilon_{kj}\leqslant \frac{C}{\lambda_j^2}.
\end{align}
Notice that
\begin{align}\label{qqqq}\lambda_i\frac{\partial\varepsilon_{ij}}{\partial\lambda_i}=-(2-\gamma)\varepsilon_{ij}\left(1-2\frac{\lambda_j}{\lambda_i}\varepsilon_{ij}^{\frac{1}{2-\gamma}}
\right).\end{align}
If $\lambda_i$ and $\lambda_j$ are comparable or $\lambda_j\leqslant\lambda_i$ (in this case, $|\lambda_j/\lambda_i|\le C$), then
\begin{align}\label{qqqq1}\lambda_i\frac{\partial\varepsilon_{ij}}{\partial\lambda_i}=-(2-\gamma)\varepsilon_{ij}(1+o(1)).\end{align}
If they are not comparable, say  $\lambda_i=o(\lambda_j)$, then
\begin{align}\label{qqqq2}\left|\lambda_i\frac{\partial\varepsilon_{ij}}{\partial\lambda_i}\right|=O(\varepsilon_{ij})\leqslant \frac{C}{\lambda_j^2}=o(\frac{1}{\lambda_i^2}).\end{align}
Thus by \eqref{case2}, there holds
\begin{align}\label{es22}
-\sum_{j\neq i} \lambda_i\frac{\partial\varepsilon_{ij}}{\partial\lambda_i}\geqslant \frac{2-\gamma}{2}\sum_{j\neq i} \varepsilon_{ij}\geqslant
\frac{2-\gamma}{2}\frac{C}{\lambda_i^2}.
\end{align}
Hence, by choosing a large $C$, it holds that
\begin{align}\label{ww}
\nonumber J'(u)\left(\lambda_i\frac{\partial w_{g_i,\lambda_i}}{\partial\lambda_i}\right)&=2\lambda(u)\left[O'(-\sum_{j\neq i}\lambda_i\frac{\partial\varepsilon_{ij}}{\partial\lambda_i})\right.\\& \nonumber\quad+\left.
\frac{2-\gamma}{4}c_2\alpha_i\frac{\Delta_{\theta_1}K(g_i)}{ K(g_i)\lambda_i^2}(1+o(1))+
o(\sum_{j\neq i} \varepsilon_{ij})\right]\\& \geqslant C_0\sum_{j\ne i}\varepsilon_{ij},
\end{align}
 where $C_0>0$ is a constant depending on $\gamma$.
Then from the first inequality of \eqref{estimate1} and \eqref{ww}, we have
\begin{align}\label{es23}
\left|J'(u)\left(\frac{1}{\lambda_i}\frac{\partial w_{g_i,\lambda_i}}{\partial g_i}\right)\right|+\frac{c'}{c}J'(u)\left(\lambda_i\frac{\partial w_{g_i,\lambda_i}}{\partial\lambda_i}\right)\geqslant c\frac{|\nabla_{\theta_1}K(g_i)|}{\lambda_i}-\frac{1}{c}\frac{1}{\lambda_i^2},\;\;\;c'=\frac{1}{C_0}.
\end{align}

In the case (3), by a simple computation,
we observe that for   $s>t$,
\begin{align}\label{0000}
-\mu_{k_s}\lambda_{k_s}\frac{\partial\varepsilon_{k_sk_t}}{\partial\lambda_{k_s}}-\mu_{k_t}\lambda_{k_t}\frac{\partial\varepsilon_{k_sk_t}}
{\partial\lambda_{k_t}}\geqslant \frac{2-\gamma}{2}\varepsilon_{k_sk_t}(1+o(1)).
\end{align}
By using \eqref{es01} and choosing a  constant $c'$ which depends on $\gamma$, there holds
\begin{align}\label{es31}
\left|J'(u)\left(\frac{1}{\lambda_i}\frac{\partial w_{g_i,\lambda_i}}{\partial g_i}\right)\right|+\frac{c'}{c}\sum_{s=1}^m\mu_{k_s}J'(u)\left(\lambda_{k_s}\frac{\partial w_{a_{k_s},\lambda_{k_s}}}{\partial\lambda_{k_s}}\right)\geqslant c\frac{|\nabla_{\theta_1} K(g_i)|}{\lambda_i}-\frac{1}{c}\frac{1}{\lambda_i^2}.
\end{align}

For a proof of \eqref{es31}, from \eqref{es01}, \eqref{qqqq} and that the functional $\lambda(u)$ has positive lower bound, we have
\begin{align*}
&\quad\sum_{k\in I_i}\mu_{k}J'(u)\left(\lambda_{k}\frac{\partial w_{g_{k},\lambda_{k}}}{\partial\lambda_{k}}\right)\\&=\sum_{k\in I_i}\mu_{k}2\lambda(u)\left[-\sum_{j\neq k}O'(\lambda_k\frac{\partial\varepsilon_{kj}}{\partial\lambda_k})+
\frac{2-\gamma}{4}c_2\alpha_k\frac{\Delta_{\theta_1} K(g_k)}{ K(g_k)\lambda_k^2}+o(\frac{1}{\lambda_k^2})+
o(\sum_{j\neq k} \varepsilon_{kj})\right]\\&\geqslant c_{\gamma}\sum_{k\in I_i}\sum_{j\neq k}\left(-\mu_{k}\lambda_k\frac{\partial\varepsilon_{kj}}{\partial\lambda_k}\right)+
\sum_{k\in I_i}\mu_{k}\lambda(u)\left(
\frac{2-\gamma}{4}c_2\alpha_k\frac{\Delta_{\theta_1} K(g_k)}{ K(g_k)\lambda_k^2}\right)\\&=\sum_{k\in I_i,j\in I_i}\sum_{j\neq k}c_{\gamma}\left(-\mu_{k}\lambda_k\frac{\partial\varepsilon_{kj}}{\partial\lambda_k}-\mu_{j}\lambda_j\frac{\partial\varepsilon_{kj}}{\partial\lambda_j}\right)+\sum_{k\in I_i,j\notin I_i}c_{\gamma}\left(-\mu_{k}\lambda_k\frac{\partial\varepsilon_{kj}}{\partial\lambda_k}\right)\\&\quad+
\sum_{k\in I_i}\mu_{k}\lambda(u)\left(
\frac{2-\gamma}{4}c_2\alpha_k\frac{\Delta_{\theta_1} K(g_k)}{ K(g_k)\lambda_k^2}\right)\\&\geqslant c_{\gamma}\sum_{k\in I_i,j\in I_i}\sum_{j\neq k}\varepsilon_{kj}+\sum_{k\in I_i}\mu_{k}\left(\sum_{\lambda_j\leqslant\lambda_k or \lambda_i \sim\lambda_j}\left(-\lambda_k\frac{\partial\varepsilon_{kj}}{\partial\lambda_k}\right)+\sum_{\lambda_k=o(\lambda_j),j\notin I_i}\left(-\lambda_k\frac{\partial\varepsilon_{kj}}{\partial\lambda_k}\right)\right)\\&\quad+
\sum_{k\in I_i}\mu_{k}\lambda(u)\left(
\frac{2-\gamma}{4}c_2\alpha_k\frac{\Delta_{\theta_1} K(g_k)}{ K(g_k)\lambda_k^2}\right)\\&\geqslant \frac{c_{\gamma}}{2}\sum_{k\in I_i}\sum_{j\neq k}\varepsilon_{kj}\geqslant \frac{c_{\gamma}}{2}\sum_{k\in I_i}\frac{C}{\lambda_k^2},
\end{align*}
where $0<c_\gamma<1$. In the last two estimates, we have used the inequalities \eqref{qqqq1} and \eqref{qqqq2} and choose the constant $C$ in \eqref{es21} large enough. $\lambda_i\sim \lambda_j$ means that $\lambda_i$ and $\lambda_j$
are comparable. It completes the the proof of claim.

\noindent{\it Proof of Lemma \ref{Wcon}.}
For the sake of simplicity, we assume $$\lambda_1\leqslant \cdots \leqslant \lambda_p.$$  We note that when we construct the vector field $W'$ satisfying the estimate \eqref{w1}, then by the same method as in \cite{bahri96} and \cite{MBA1996}, we can prove \eqref{w2}. So in the following we only need to construct the vector field $W'$ satisfies the inequality  \eqref{w1}.
We divide it into four cases.

\noindent{\it Case 1.} Suppose there holds
\begin{align}\label{ncase1}
\frac{|\nabla_{\theta_1}K(g_1)|}{\lambda_1}>\frac{2}{c^2\lambda_1^2}.
\end{align}
In this case, by \eqref{esaa} we have
\begin{align}\label{nes11}
\left|J'(u)\left(\frac{1}{\lambda_1}\frac{\partial w_{g_1,\lambda_1}}{\partial g_1}\right)\right|+\frac{c'}{c}\sum_{k\in I_1}\mu_kJ'(u)\left(\lambda_k\frac{\partial w_{g_k,\lambda_k}}{\partial\lambda_k}\right)\geqslant \frac{c}{4}\frac{|\nabla_{\theta_1} K(g_1)|}{\lambda_1}+\frac{1}{4c}\frac{1}{\lambda_1^2}.
\end{align}
Combined with (\ref{es31}), (\ref{nes11}) and (\ref{es01}), for any $i$,  we reach
\begin{align}\label{nes12}
\Gamma_i:&=J'(u)\left(\lambda_i\frac{\partial w_{g_i,\lambda_i}}{\partial \lambda_i}\right)+\left|J'(u)\left(\frac{1}{\lambda_i}\frac{\partial w_{g_i,\lambda_i}}{\partial g_i}\right)\right|+\frac{\tilde{c}}{c}\sum_{k\in I_i}\mu_kJ'(u)\left(\lambda_k\frac{\partial w_{g_k,\lambda_k}}{\partial\lambda_k}\right)\nonumber \\&\quad+4\left(\left|J'(u)\left(\frac{1}{\lambda_1}\frac{\partial w_{g_1,\lambda_1}}{\partial g_1}\right)\right|+\frac{\tilde{c}}{c}\sum_{k\in I_1}\mu_kJ'(u)\left(\lambda_k\frac{\partial w_{g_k,\lambda_k}}{\partial\lambda_k}\right)\right)\nonumber\\&\geqslant -\sum_{j\neq i}O'(\lambda_i\frac{\partial\varepsilon_{ij}}{\partial\lambda_i})+
\frac{B}{\lambda_i^2}(1+o(1))+
o(\sum_{j\neq i} \varepsilon_{ij})+c\frac{|\nabla_{\theta_1}K(g_i)|}{\lambda_i}-\frac{1}{c}\frac{1}{\lambda_i^2}\nonumber \\&\quad+
4\left(\frac{c}{4}\frac{|\nabla_{\theta_1}K(g_1)|}{\lambda_1}+\frac{1}{4c}\frac{1}{\lambda_1^2}\right)\nonumber \\&\geqslant-\sum_{j\neq i}O'(\lambda_i\frac{\partial\varepsilon_{ij}}{\partial\lambda_i})+
\frac{B}{\lambda_i^2}+c\frac{|\nabla_{\theta_1}K(g_i)|}{\lambda_i}+o(\sum_{j\neq i} \varepsilon_{ij}).
\end{align}

Here and in sequel we denote $B>0$  a  constant which may vary in different places,
We define $\nu_k=2^{k-1}$. By simple computation, we have
\begin{align*}
-\left(\nu_j\lambda_j\frac{\partial\varepsilon_{ij}}{\partial\lambda_j}+\nu_i\lambda_i\frac{\partial\varepsilon_{ij}}{\partial\lambda_i}\right)\geqslant D_\gamma\varepsilon_{ij},\
\textmd{if}\ j>i, \;D_\gamma=\frac{2-\gamma}2.
\end{align*}
So we get
\begin{align}\label{nes13}
\sum_{i=1}^p\nu_i\Gamma_i\geqslant B\sum_{i\neq j}\varepsilon_{ij}+
\sum_{i=1}^p\frac{B}{\lambda_i^2}+B\sum_{i=1}^p\frac{|\nabla_{\theta_1}K(g_i)|}{\lambda_i}.
\end{align}
This can be rewritten in the following form
\begin{align}\label{nes14}
J'(u)\left(\sum_{i=1}^p\gamma_i\lambda_i\frac{\partial w_{g_i,\lambda_i}}{\partial \lambda_i}\right)+\sum_{i=1}^p\beta_i\left|J'(u)\left(\frac{1}{\lambda_i}\frac{\partial w_{g_i,\lambda_i}}{\partial g_i}\right)\right|\geqslant B\left(\sum_{i\neq j}\varepsilon_{ij}+
\sum_{i=1}^p\frac{1}{\lambda_i^2}+\sum_{i=1}^p\frac{|\nabla_{\theta_1}K(g_i)|}{\lambda_i}\right),
\end{align}
where $\gamma_i, \beta_i$ are bounded nonnegative constants depending on $\mu_i,\;\nu_i$ and $\gamma$.

We now define the vector field by
\begin{align}\label{nes15}
W'=-\left(\sum_{i=1}^p\gamma_i\lambda_i\frac{\partial w_{g_i,\lambda_i}}{\partial \lambda_i}\right)-\sum_{i=1}^p\beta_i \left(J'(u)\left(\frac{1}{\lambda_i}\frac{\partial w_{g_i,\lambda_i}}{\partial g_i}\right)\right)\cdot\frac{\nabla_{\theta_1}K(g_i)}{|\nabla_{\theta_1}K(g_i)|}.
\end{align}
From the estimate \eqref{nes14}, we have \eqref{w1}.

\noindent{\it Case 2.} Suppose there holds
\begin{align}\label{ncase2}
-\sum_{j\neq 1}O'(\lambda_1\frac{\partial\varepsilon_{1j}}{\partial\lambda_1})>\frac{4(|\Delta_{\theta_1} K(g_1)|+1)}{\lambda_1^2}.
\end{align}
In this case, as in \eqref{nes12} we define
\begin{align}\label{nes21}
\Gamma_i=&J'(u)\left(\lambda_i\frac{\partial w_{g_i,\lambda_i}}{\partial \lambda_i}\right)+\left|J'(u)\left(\frac{1}{\lambda_i}\frac{\partial w_{g_i,\lambda_i}}{\partial g_i}\right)\right|+\frac{\tilde{c}}{c}\sum_{k\in I_i}\mu_kJ'(u)\left(\lambda_k\frac{\partial w_{g_k,\lambda_k}}{\partial\lambda_k}\right)\nonumber \\&+BJ'(u)\left(\lambda_1\frac{\partial w_{g_1,\lambda_1}}{\partial \lambda_1}\right).
\end{align}
A similar construct of $W'$ can be done and the proof of \eqref{w1} is repeated as the case 1 word by word by some mirror modifications.

\noindent{\it Case 3.} Suppose there holds
\begin{align}\label{ncase3}
\sum_{j\neq i}\varepsilon_{ij}\geqslant\frac{C}{\lambda_1^2},\  C \ \textmd{a suitable constant}.
\end{align}
In this case, we define
\begin{align}\label{nes31}
\Gamma_i=J'(u)\left(\lambda_i\frac{\partial w_{g_i,\lambda_i}}{\partial \lambda_i}\right)+\left|J'(u)\left(\frac{1}{\lambda_i}\frac{\partial w_{g_i,\lambda_i}}{\partial g_i}\right)\right|+\frac{\tilde{c}}{c}\sum_{k\in I_i}\mu_kJ'(u)\left(\lambda_k\frac{\partial w_{g_k,\lambda_k}}{\partial\lambda_k}\right).
\end{align}
Since there holds
\begin{align*}
-\sum_{i=1}^p\nu_i\sum_{j\neq i}O'(\lambda_i\frac{\partial\varepsilon_{ij}}{\partial\lambda_i})\geqslant B\sum_{j\neq i}\varepsilon_{ij},
\end{align*}
so we can define vector field $W'$ and give a similar proof of \eqref{w1} as in case 1.

\noindent{\it Case 4.} Suppose  there holds
\begin{align}\label{ncase4}
\sum_{j\neq i}\varepsilon_{ij}<\frac{C}{\lambda_1^2}, \ \ \
\frac{|\nabla_{\theta_1}K(g_1)|}{\lambda_1}<\frac{2}{c^2\lambda_1^2}.
\end{align}
the above proof extends as follows.

\noindent{\it Subcase 1.} Suppose in the sequence $\lambda_1\leqslant \cdots \leqslant \lambda_p$, there exists $i_1$ such that for some   $0\le r\le p-i_1$, there holds

\begin{align}\label{subcase}
\sum_{s=0}^r\sum_{j\leqslant i_1+r,j\neq i_1+s}\varepsilon_{i_1+r,j}\geqslant\frac{C}{\lambda_{i_1}^2}, \ \textmd{or}\ \
\frac{|\nabla_{\theta_1}K(g_{i_1})|}{\lambda_{i_1}}\geqslant\frac{2}{c^2\lambda_{i_1}^2}.
\end{align}
We note that for a choice of $(i_1,r_0)$ satisfying (\ref{subcase}), then all of $(i_1,r)$ with $r_0\le r\le p-i_1$ satisfies (\ref{subcase}).
Similarly to the case 1 and case 3, we can define  a vector field $W(i_1,r)$ in $\textmd{span}_{i=i_1}^{i_1+r}\{\frac{\partial w_{g_i,\lambda_i}}{\partial \lambda_i},\frac{\partial w_{g_i,\lambda_i}}{\partial g_i}\}$ such that
$$\|W(i_1,r)\|\leqslant C$$
and
\begin{align}\label{subes}
-J'(u)W(i_1,r)&\geqslant B\left(\sum_{s=0}^r\sum_{j\leqslant i_1+r,j\neq i_1+s}\varepsilon_{i_1+s,j}+
\sum_{s=0}^r\frac{1}{\lambda_{i_1+s}^2}\right.\nonumber \\&\quad+\left.\sum_{s=0}^r\frac{|\nabla_{\theta_1} K(g_{i_1+s})|}{\lambda_{i_1+s}}-\frac{1}{\bar{c}}\sum_{s=0}^r\sum_{j\geqslant i_1+r+1,}\varepsilon_{i_1+s,j}\right).
\end{align}

Assume  $i_1$ is the smallest subscript satisfying (\ref{subcase}).
Then by choosing $r=p-i_1$, we have
\begin{align}\label{ecase1es1}
-J'(u)W(i_1,p-i_1)\geqslant B\left(\sum_{k=i_1,j\neq k}^p\varepsilon_{jk}+
\sum_{j\geqslant i_1}\frac{1}{\lambda_{j}^2}+\sum_{j\geqslant i_1}\frac{|\nabla_{\theta_1}K(g_{j})|}{\lambda_{j}}\right).
\end{align}

If $i_1=1$, we obtain the result of \eqref{w1}.

Otherwise, for integer $l\in [1,i_1)$, there holds
\begin{align}\label{ecase1sub2}
\sum_{k\geqslant l}\sum_{j\neq k}\varepsilon_{jk}\leqslant\frac{C}{\lambda_{l}^2} \ \textmd{and}\ \
\frac{|\nabla_{\theta_1}K(g_{l})|}{\lambda_l}\leqslant\frac{2}{c^2\lambda_l^2}.
\end{align}

From (\ref{ecase1es1}) and (\ref{ecase1sub2}), then \eqref{w1} is true if
\begin{align}\label{ecase1sub21}
\sum_{k=i_1,j\neq k}^p\varepsilon_{jk}+
\sum_{j\geqslant i_1}\frac{1}{\lambda_{j}^2}+\sum_{j\geqslant i_1}\frac{|\nabla_{\theta_1}K(g_{j})|}{\lambda_{j}}\ge  \frac{B}{\lambda_{1}^2} \ \textmd{for some}\ B>0.
\end{align}

If (\ref{ecase1sub21}) is not hold, then there holds
\begin{align}\label{ecase1sub211}
\sum_{k=i_1,j\neq k}^p\varepsilon_{jk}+
\sum_{j\geqslant i_1}\frac{1}{\lambda_{j}^2}+\sum_{j\geqslant i_1}\frac{|\nabla_{\theta_1}K(g_{j})|}{\lambda_{j}}= o\left(\frac{1}{\lambda_{1}^2}\right).
\end{align}

Combining with (\ref{ecase1sub2}), we have for $j\leqslant i_1-1$,
\begin{align*}
\lambda_j|\nabla_{\theta_1}K(g_{j})|\leqslant\frac{2}{c^2}, \ \ |\nabla_{\theta_1}K(g_{j})|=o(1).
\end{align*}
These imply that: for $j\leqslant i_1-1, g_j$ is close to a critical point of $K$ which we denoted by $\eta_j$,
$\lambda_j\textmd{d}(g_j,\eta_j)=O(1).$

If $i\leqslant i_1-1, j\leqslant i_1-1, \eta_i=\eta_j$, we have
$$|\textmd{inf}(\lambda_i, \lambda_j)\textmd{d}(g_i,g_j)|=O(1).$$
So if $j< i,$
\begin{align*}
o(1)=\varepsilon_{ij}\geqslant C(\frac{\lambda_j}{\lambda_i})^{2-\gamma} \ \textmd{and}\ \varepsilon_{ij}\leqslant\frac{C}{\lambda_i^2}=o(\frac{1}{\lambda_j^2})
\end{align*}
 and for $1<j\le i_1-1$ the holds $\frac{1}{\lambda_j^2}=o(\frac{1}{\lambda_1^2})$.
If $i\leqslant i_1-1, j\leqslant i_1-1, \eta_i\neq \eta_j$, we have
$$\varepsilon_{ij}=O(\frac{1}{\lambda_i\lambda_j})^{2-\gamma}=o(\frac{1}{\lambda_1^2}).$$
Thus together with \eqref{ecase1sub211}, we have
$$\sum_{i\neq j}\varepsilon_{ij}=o(\frac{1}{\lambda_1^2}).$$

Combining \eqref{ecase1sub2} and \eqref{ecase1sub211}, we have
\begin{align*}
\sum_{i=1}^p\frac{|\nabla_{\theta_1}K(g_i)|}{\lambda_i}&=\sum_{i\geqslant i_1}\frac{|\nabla_{\theta_1}K(g_i)|}{\lambda_i}+\sum_{i< i_1}\frac{|\nabla_{\theta_1} K(g_i)|}{\lambda_i}\leqslant
\frac{c}{\lambda_1^2}.
\end{align*}

So
\begin{align*}
\frac{1}{\lambda_1^2}\geqslant B\left(\sum_{i=1}^p\frac{|\nabla_{\theta_1}K(g_i)|}{\lambda_i}+\sum_{i=1}^p\frac{1}{\lambda_i^2}+\sum_{i\neq j}\varepsilon_{ij}\right).
\end{align*}

Since $g_1$ is close to a critical point of $K$ which we denoted by $\eta_1$ and $\Delta_{\theta_1} K(\eta_1)\ne 0$, we get
\begin{align*}
-J'(u)\left(\lambda_1\frac{\partial w_{g_1,\lambda_1}}{\partial\lambda_1}\right)&=O'(\sum_{j\neq 1}\lambda_1\frac{\partial\varepsilon_{1j}}{\partial\lambda_1})-
c\frac{\Delta_{\theta_1} K(\eta_1)}{\lambda_1^2}(1+o(1))\\&=-
c\frac{\Delta_{\theta_1} K(\eta_1)}{\lambda_1^2}+o(\frac{1}{\lambda_1^2}).
\end{align*}

If $-\Delta_{\theta_1} K(\eta_1)\geqslant c>0$, it holds that
\begin{align*}
-J'(u)\left(\lambda_1\frac{\partial w_{g_1,\lambda_1}}{\partial\lambda_1}\right)&=-
c\frac{\Delta_{\theta_1} K(\eta_1)}{\lambda_1^2}+o(\frac{1}{\lambda_1^2})\\&\geqslant c'\left(\sum_{i=1}^p\frac{|\nabla_{\theta_1} K(g_i)|}{\lambda_i}+\sum_{i=1}^p\frac{1}{\lambda_i^2}+\sum_{i\neq j}\varepsilon_{ij}\right).
\end{align*}
Now we define the vector field
\begin{align}\label{0-0}
W'=\lambda_1\frac{\partial w_{g_1,\lambda_1}}{\partial \lambda_1},
\end{align}
which satisfies \eqref{w1}.

If $-\Delta_{\theta_1} K(\eta_1)\leqslant -c<0$, we define
\begin{align}\label{1-0}
W'=-\lambda_1\frac{\partial w_{g_1,\lambda_1}}{\partial \lambda_1},
\end{align}
which also satisfies \eqref{w1}.

\noindent{\it Subcase 2.} Assume that indices $i_1$ satisfying (\ref{subcase}) do not exist, i.e., for any $l\in\{1, \cdots, p\}$
\begin{align}\label{ecase2}
\sum_{k=l}^p\sum_{i\neq k}\varepsilon_{ik}\leqslant\frac{C}{\lambda_{l}^2} \ \textmd{and}\ \
\frac{|\nabla_{\theta_1}K(g_{l})|}{\lambda_{l}}\leqslant\frac{2}{c^2\lambda_{l}^2}.
\end{align}

By a direct argument, when $i<j, \eta_i=\eta_j$, we get
$$\textmd{inf}(\lambda_i, \lambda_j)\textmd{d}(g_i,g_j)=O(1).$$
Thus under the condition \eqref{ecase2}, for some $i< j$
\begin{align*}
o(1)=\varepsilon_{ij}\geqslant C(\frac{\lambda_i}{\lambda_j})^{2-\gamma} \ \textmd{and}\ \varepsilon_{ij}\leqslant\frac{C}{\lambda_j^2}=o(\frac{1}{\lambda_i^2}).
\end{align*}
The construct of vector field is same as \eqref{0-0} or \eqref{1-0} in the previous subcase.
In fact, we have reach that if two $\lambda$'s for example $\lambda_i$ and $\lambda_j$ are not comparable, the vector field $W'$ can be defined to satisfy \eqref{w1}.
So in this subcase,  we can assume that $\textmd{inf}_{i\neq j}\textmd{d}(g_i,g_j)\geqslant d_0>0$ and all the $\lambda$'s are comparable. Thus we have $$\sum_{j\neq i}\lambda_i\frac{\partial\varepsilon_{ij}}{\partial\lambda_i}=\sum_{j\neq i}-(2-\gamma)\varepsilon_{ij}(1+o(1))=o(\frac{1}{\lambda_1^2}).$$
Therefore, we have
\begin{align*}
-J'(u)\left(\lambda_i\frac{\partial w_{g_i,\lambda_i}}{\partial\lambda_i}\right)=-
c\frac{\Delta_{\theta_1} K(\eta_i)}{\lambda_i^2}+o(\frac{1}{\lambda_i^2}).
\end{align*}
If for some $\eta_i$ satisfies $-\Delta_{\theta_i} K(\eta_i)\leqslant -c<0$, we define
\begin{align}\label{1-0}
W'=-\lambda_i\frac{\partial w_{g_i,\lambda_1}}{\partial \lambda_i}.
\end{align}
If for all $\eta_i$, $-\Delta_{\theta_i} K(\eta_i)\geqslant c>0$,
Now the construct of vector field is same as $$W'=\sum_{i=1}^p\lambda_i\frac{\partial w_{g_i,\lambda_i}}{\partial\lambda_i}.$$

Since in the cases 1-4, we can adjust some constants to insure that the  union of these four case is the whole discussed space, by using a partition of unity, we can define the final vector field $W'$ satisfying \eqref{w1} at all.
 The proof of Lemma \ref{Wcon} is complete.

 \hfill $\Box$

\begin{lemma}\label{4.3}
For any $u=\sum_{i=1}^p\alpha_iw_{g_i,\lambda_i}\in V(p,\varepsilon') (\varepsilon'<\frac{\varepsilon}{2})$, there exist $\tilde g=(\tilde g_1,\cdots,\tilde g_p)$ and $\tilde \lambda=(\tilde \lambda_1,\cdots,\tilde \lambda_p)$ such that
\begin{align}\label{1-00}J\left(\sum_{i=1}^p\alpha_iw_{g_i,\lambda_i}+\bar{v}(\alpha,g,\lambda)\right)=J\left(\sum_{i=1}^p\alpha_i w_{\tilde{g_i},\tilde{\lambda_i}}\right)\end{align}
and the following two statements
\begin{align}\label{1-01}\sum_{i\neq j} \tilde{\varepsilon}_{ij}+\sum_{i=1}^p\frac{1}{\tilde{\lambda}_i^2}\rightarrow 0\Longleftrightarrow\sum_{i\neq j} \varepsilon_{ij}+\sum_{i=1}^p\frac{1}{\lambda_i^2}\rightarrow 0,\end{align}
\begin{align}\label{1-02}\textmd{d}(\tilde{g}_i,g_i)\rightarrow 0 \ \textmd{as} \ \sum_{i\neq j} \varepsilon_{ij}+\sum_{i=1}^p\frac{1}{\lambda_i^2}\rightarrow 0.\end{align}
\end{lemma}
\pf
The proof is similar to the one given in \cite{MBA1996} and \cite{3CR02}.

By lemma \ref{Wcon}, the vector field $W'$ is Lipschitz. Hence, there is a 1-parameter group $h_s$ generated by $W'$ satisfying
\begin{align*}
\begin{cases}
\frac{\partial}{\partial s}h_s(\sum_{i=1}^p\alpha_iw_{g_i,\lambda_i})=W'\left(h_s(\sum_{i=1}^p\alpha_iw_{g_i,\lambda_i})\right),\\
h_0(\sum_{i=1}^p\alpha_iw_{g_i,\lambda_i})=\sum_{i=1}^p\alpha_iw_{g_i,\lambda_i}.
\end{cases}
\end{align*}

For different critical points $\eta_i'$s of $K$ with $-\Delta_{\theta_1} K(\eta_i)\geqslant c>0$, $i\in\{1,\cdots,p\}$ and $\delta<\frac 12\min\{dist(\eta_i,\eta_j)\}$, we define $V_\delta(\eta_1,\cdots,\eta_p)$ to be the set of $(g,\lambda)$ satisfying  $g_i\in B_{\delta}(\eta_i)$,  $i=1,\cdots,p$.

$J\left(h_s(\sum_{i=1}^p\alpha_iw_{g_i,\lambda_i})\right)$ and $J\left(h_s(\sum_{i=1}^p\alpha_iw_{g_i,\lambda_i}+\bar v(s))\right)$ are decreasing functions of $s$. Since $J(\sum_{i=1}^p\alpha_iw_{g_i,\lambda_i}+\bar v)\leqslant J(\sum_{i=1}^p\alpha_iw_{g_i,\lambda_i})$, there is at most one solution of the equation
\begin{align}\label{916}
J\left(h_s\left(\sum_{i=1}^p\alpha_iw_{g_i,\lambda_i}\right)\right)=J\left(\sum_{i=1}^p\alpha_iw_{g_i,\lambda_i}+\bar v\right).
\end{align}
By using Lemma \ref{Wcon} and a similar proof as in \cite{bahri1991}, we can see that the flow line $h_s(\sum_{i=1}^p\alpha_iw_{g_i,\lambda_i})$ satisfies the (PS) condition  if $u_0=\sum_{i=1}^p\alpha_iw_{g_i,\lambda_i}\notin V_\delta(\eta_1,\cdots,\eta_p)$. {\it i.e.}, for a small fixed $\varepsilon_0>0$, there is $\varepsilon_1>0$ such that $h_s(\sum_{i=1}^p\alpha_iw_{g_i,\lambda_i})$ remains outside $V(p,\varepsilon_1)$ when $s\ge \varepsilon_0$.

The cases in which there could be no solution of (\ref{916}) are $h_s\left(\sum_{i=1}^p\alpha_iw_{g_i,\lambda_i}\right)$ exits from $V(p,\varepsilon_1)$ or the decreasing flow goes to critical points at infinity.

If $h_s\left(\sum_{i=1}^p\alpha_iw_{g_i,\lambda_i}\right)$ exits from $V(p,\varepsilon_1)$, the flow line have to travel from $V(p,\frac{\varepsilon_1}{2})$ to $V(p,\varepsilon_1)$. By Lemma \ref{Wcon}, $\frac{\partial}{\partial s}J\left(h_s\left(\sum_{i=1}^p\alpha_iw_{g_i,\lambda_i}\right)\right)$ is lower bounded by a constant $\delta_1>0$ and $d(\partial V(p,\varepsilon_1),\partial V(p,\frac{\varepsilon_1}{2}))=\delta_2>0$. Since $\|W'\|\leqslant C$, then $J\left(h_s\left(\sum_{i=1}^p\alpha_iw_{g_i,\lambda_i}\right)\right)$ decreases at least $\frac{\delta_1\delta_2}{C}$. However,
$$J\left(\sum_{i=1}^p\alpha_iw_{g_i,\lambda_i}\right)-J\left(\sum_{i=1}^p\alpha_iw_{g_i,\lambda_i}+\bar v\right)\rightarrow 0, \;\;\; \varepsilon\rightarrow 0.$$
We can choose $\varepsilon>0$ sufficiently small, there is a solution of (\ref{916}).

If $\sum_{i=1}^p\alpha_iw_{g_i,\lambda_i}\in V_\varepsilon(\eta_1,\cdots,\eta_p)$, by Lemma \ref{Wcon}, it will take an infinity time for the flow to go to infinity. Therefore, at least for a subsequence ${s_k}$,$s_k\rightarrow +\infty$,
$$\varepsilon_{ij}(s_k)+\sum_{i=1}^p\frac{1}{\lambda_i^2}(s_k)\rightarrow 0.$$
This implies $\bar v(s_k)\rightarrow 0$ and $J\left(u(s_k)\right)-J\left(\bar u(s_k)\right)\rightarrow 0$. Thus
\begin{align*}\textmd{liminf}_{s\rightarrow +\infty}J(u(s))=\textmd{liminf}_{s\rightarrow +\infty}J(\bar u(s))<J(\bar u).
\end{align*}
By continuity, (\ref{916}) must have a solution.

Similarly,  we consider the vector field $-W'$ and the flow line $h_{-s}(\sum_{i=1}^p\alpha_i w_{\tilde{g_i},\tilde{\lambda_i}})$.
It is easy to know that  there is a unique solution for
$$J\left(h_{-s}(\sum_{i=1}^p\alpha_i w_{\tilde{g_i},\tilde{\lambda_i}})+\bar{v}(h_{-s}(\sum_{i=1}^p\alpha_i w_{\tilde{g_i},\tilde{\lambda_i}}
))\right)
=J\left(\sum_{i=1}^p\alpha_i w_{\tilde{g_i},\tilde{\lambda_i}}\right).$$

Set
$$h_s\left(\sum_{i=1}^p\alpha_iw_{g_i,\lambda_i}\right)=\sum_{i=1}^p\alpha_iw_{g_i(s),\lambda_i(s)},$$ and take $(\tilde{g_i},\tilde{\lambda_i})=(g_i(s),\lambda_i(s))$, we have \eqref{1-00}.

As for \eqref{1-01} and \eqref{1-02}, we note that
$$W'=\sum_{i=1}^p\alpha_i\left(\frac{1}{\lambda_i(s)}\frac{\partial w_{g_i(s),\lambda_i(s)}}{\partial g_i(s)}\right)\left(\lambda_i(s)\dot{g}_i(s)\right)+\sum_{i=1}^p\alpha_i\left(\lambda_i(s)\frac{\partial w_{g_i(s),\lambda_i(s)}}{\partial \lambda_i(s)}\right)\frac{\dot{\lambda}_i(s)}{\lambda_i(s)},$$
where $\dot{g}_i(s)$ and $\dot{\lambda}_i(s)$ denote the action of $W'$ on the variables $g_i$ and $\lambda_i$. We have $|\lambda_i\dot{g}_i|\leqslant C$,$|\frac{\dot{\lambda}_i}{\lambda_i}|\leqslant C, i=1,\cdots,p$. Then
\begin{align*}
\left|\frac{\partial\varepsilon_{ij}(s)}{\partial s}\right|=\left|\frac{\partial\varepsilon_{ij}}{\partial \lambda_i}\frac{\partial\lambda_{i}(s)}{\partial s}+\frac{\partial\varepsilon_{ij}}{\partial \lambda_j}\frac{\partial\lambda_{j}(s)}{\partial s}+\frac{\partial\varepsilon_{ij}}{\partial g_i}\frac{\partial g_{i}(s)}{\partial s}+\frac{\partial\varepsilon_{ij}}{\partial g_j}\frac{\partial g_{j}(s)}{\partial s}\right|\leqslant C\varepsilon_{ij}(s).
\end{align*}
Thus, $$e^{-Cs}\varepsilon_{ij}\leqslant \varepsilon_{ij}(s)\leqslant e^{Cs}\varepsilon_{ij},\\
e^{-Cs}\leqslant \frac{\lambda_{i}(s)}{\lambda_{i}(0)} \leqslant e^{Cs},\\
|g_{i}(s)-g_i|\leqslant \frac{e^{Cs}}{\lambda_{i}(0)}.
$$
Since the $s$ satisfying (\ref{916}) is bounded, we get  \eqref{1-01} and \eqref{1-02}.\hfill $\Box$\\\

Following from Lemma \ref{4.3}, for any $\varepsilon_1>0$ small, there are $\varepsilon>0$ and $\varepsilon_2>0$ such that
$$V(p,\varepsilon)\xrightarrow{{h_s}}V(p,\varepsilon_1)\xrightarrow{h_{-s}}V(p,\varepsilon_2)\supseteq V(p,\varepsilon_1).$$
From Lemma \ref{4.1}, Lemma \ref{4.3} and this fact, we can proof Theorem \ref{covering}.

\section{The proof of main theorem}
For technical reasons, we introduce for $\varepsilon_0>0$ small enough, the following subset of $\Sigma$
as $V_{\varepsilon_0}(\Sigma^+)=\{u\in\Sigma, |u^-|_{L^{\frac{4}{2-\gamma}}}<\varepsilon_0\}$.

By Theorem \ref{covering}, there is a covering $\{O_l\}$ of the base space of the bundle $V(p,\varepsilon)$ such that  Theorem \ref{covering} holds on each $\{O_l\}$. For $u=\sum_{i=1}^p\alpha_i w_{g_i,\lambda_i}+v \in V(p,\varepsilon)$, we consider transformation of coordinates $\varphi:(\alpha,g,\lambda,v)\rightarrow (\alpha,\tilde{g},\tilde{\lambda},V_l)$, so that \eqref{zz} holds.

We first define a vector field on $V(p,\varepsilon)$ by using a partition of unity $\eta_l$ on the base space of $(\alpha, \tilde g, \tilde \lambda, V)$ as
 $X=W-\sum_{l=1}^m\eta_l V_l$, where $W(\alpha,\tilde{g},\tilde{\lambda},V_l)=W'(\alpha,\tilde{g},\tilde{\lambda})$. Then the vector field
  in the variables $(\alpha,g,\lambda,v)$ is defined as $Z=X\circ \varphi$. By direct computation, in every open set $O_l$
  we have $$-J'\left(\sum_{i=1}^p
\alpha_i w_{g_i,\lambda_i}+v\right)(Z)=-J'\left(\sum_{i=1}^p
\alpha_i w_{\tilde g_i,\tilde\lambda_i}\right)(W')+J''\left(u_0\right)V_l\cdot V_l+O(\|v\|^2),$$
where $u_0=\sum_{i=1}^p\alpha_i w_{g_i,\lambda_i}+\bar{v}(\alpha,g,\lambda)$.

We remind that $V_l$ orthogonal to $w_{g_i,\lambda_i},\frac{\partial w_{g_i,\lambda_i}}{\partial \lambda_i},\frac{\partial w_{g_i,\lambda_i}}{\partial g_i}$, so by computation,
\begin{align*}
J''\left(u_0\right)V_l\cdot V_l&=2\lambda(u_0)||V_l||^2-4\lambda(u_0)^{\frac{4-\gamma}{2-\gamma}}\int_{\mathbb{S}^{3}}P_{\gamma}u_0V_l  \int_{\mathbb{S}^{3}}Ku_0^{\frac{2+\gamma}{2-\gamma}}V_l\\&\quad+\frac{4(4-\gamma)}{2-\gamma}\lambda(u_0)^{\frac{6-\gamma}{2-\gamma}} \left(\int_{\mathbb{S}^{3}}Ku_0^{\frac{2+\gamma}{2-\gamma}}V_l\right)^2-\frac{2(2+\gamma)}{2-\gamma}\lambda(u_0)^{\frac{4-\gamma}{2-\gamma}} \int_{\mathbb{S}^{3}}Ku_0^{\frac{2\gamma}{2-\gamma}}V_l^2.
\end{align*}
Using the estimates in Appendix A and $Q(v,v)$ is positive definite, we obtain
$J''\left(u_0\right)V_l\cdot V_l$ is positive definite. So in $V(p,\varepsilon)$, if $\varepsilon$ small enough, there holds
$$-J'\left(\sum_{i=1}^p
\alpha_i w_{g_i,\lambda_i}+v\right)(Z)\ge C\left(\sum_{i=1}^p(\frac{|\nabla_{\theta_1}K(\tilde{g}_i)|}{\tilde{\lambda}_i}+\frac{1}{\tilde{\lambda}_i^2})+\sum_{i\neq j} \tilde{\varepsilon}_{ij}\right).$$

We now suppose that the functional $J$ has no critical point and  there holds $\|J'(u)\|\ge C_\varepsilon>0$ for $u\notin V(p,\frac{\varepsilon}{2})$. On $V(p,\frac{\varepsilon}{2})$, define the vector field $-J'$, and then also via a partition of unity of the two sets $V(p,\varepsilon)$ and $V(p,\frac{\varepsilon}{2})$ to
 build a global vector field $\bar Z(u)$ on $V_{\varepsilon_0}(\Sigma^+)$ from $Z$ and $-J'$. It is easy to see
$$-J'(u)(\bar Z)\geqslant C\left(\sum_{i=1}^p(\frac{|\nabla_{\theta_1}K(\tilde{g}_i)|}{\tilde{\lambda}_i}+\frac{1}{\tilde{\lambda}_i^2})+\sum_{i\neq j} \tilde{\varepsilon}_{ij}\right).$$

It is important to insure that any flow line generated by the vector field $\bar Z$ with initial
condition $u\in V_{\varepsilon_0}(\Sigma^+)$ remains in $V_{\varepsilon_0}(\Sigma^+)$. We have the following lemma.
\begin{lemma}
$V_{\varepsilon_0}(\Sigma^+)$ is invariant under the flow generated by $\bar Z(u)$ .
\end{lemma}
\pf
It is sufficient to prove that $V_{\varepsilon_0}(\Sigma^+)$ is invariant under the negative gradient flow
of $J$.

Suppose $u_0\in V_{\varepsilon_0}(\Sigma^+)$,
\begin{align*}
\begin{cases}
\frac{\textmd{d}u(s)}{\textmd{d}s}=-2\lambda(u(s)) u(s)+2\lambda(u(s))^{\frac{4-\gamma}{2-\gamma}}P_\gamma^{-1}Ku(s)^{\frac{2+\gamma}{2-\gamma}},\\
u(0)=u_0.
\end{cases}
\end{align*}
Then
$$e^{\int_0^s2\lambda(u(t))dt}u(s)=u_0+\int_0^se^{\int_0^t2\lambda(u(\zeta))d\zeta}2\lambda(u(t))^{\frac{4-\gamma}{2-\gamma}}P_\gamma^{-1}Ku(t)^{\frac{2+\gamma}{2-\gamma}}\textmd{d}t.$$
Therefore,
\begin{align*}
u(s)&=e^{-\int_0^s2\lambda(u(t))dt}u_0^++e^{-\int_0^s2\lambda(u(t))dt}\int_0^se^{\int_0^t2\lambda(u(\zeta))d\zeta}2\lambda(u(t))^{\frac{4-\gamma}{2-\gamma}}P_\gamma^{-1}Ku(t)^{\frac{2+\gamma}{2-\gamma}}
\textmd{d}t\\&\quad-e^{-\int_0^s2\lambda(u(t))dt}u_0^-.
\end{align*}
Hence,
\begin{align*}
u^-(s)\leqslant e^{-\int_0^s2\lambda(u(t))dt}u_0^-.
\end{align*}
Set
\begin{align*}
f(s)=e^{-\frac{4}{2-\gamma}\int_0^s2\lambda(u(t))dt}|u_0^-|_{L^{\frac{4}{2-\gamma}}}^{\frac{4}{2-\gamma}}.
\end{align*}
Then $|u^-(s)|_{L^{\frac{4}{2-\gamma}}}^{\frac{4}{2-\gamma}}\leqslant f(s)$,
\begin{align*}
f'(s)=-\frac{8\lambda}{2-\gamma} e^{-\frac{4}{2-\gamma}\int_0^s2\lambda(u(t))dt}|u_0^-|_{L^{\frac{4}{2-\gamma}}}^{\frac{4}{2-\gamma}}\leqslant 0.
\end{align*}
Therefore,
$$ |u(s)^-|_{L^{\frac{4}{2-\gamma}}}<\varepsilon_0  \ \textmd{for all}\ s> 0.$$\hfill $\Box$\\

Next, we study the concentration phenomenon of the functional $J$.
\begin{lemma}
Assume that (\ref{main}) has no solution. Then, the set of critical point at
infinity of $J$ in $\Sigma^+$ lie in $\cup_{p=1}^mV_\delta(\eta_1,\cdots,\eta_p)$.
\end{lemma}
\pf
From the fact that outside $\cup_p  V(p,\frac{3}{4}\varepsilon)$, $-J'(u)(\bar Z)\geqslant C$ and lemma \ref{PSfail},
 we know that if $u_0\in V_{\varepsilon_0}(\Sigma^+)$,
there is $p\in N^*$  and $s_0>0$ such that if $\eta(s, u_0)$ denotes the flow line of the vector field
$Z$ with initial condition $u_0$, that is $\eta(s, u_0)$ satisfies
\begin{align*}
\begin{cases}
\frac{\partial}{\partial s}\eta(s, u_0)=\bar Z\left(\eta(s, u_0)\right),\\
\eta(0, u_0)=u_0,
\end{cases}
\end{align*}
$\eta(s, u_0)$ is in $ V(p,\frac{3}{4}\varepsilon)$ for $s\geqslant s_0$, since in $\in V_{\varepsilon_0}(\Sigma^+)$, $J(u)>0$.

Assume that for any $s\geqslant s_0$, $\eta(s, u_0)$ is in $ V(p,\frac{3}{4}\varepsilon)$ but outside $V_\delta(\eta_1,\cdots,\eta_p)$. It means that $\eta(s, u_0)$ satisfies (PS) condition for $s\geqslant s_0$.
By the construction of the vector field $W'$,  in
 all the cases of the vector field defined with negative coefficients in direction $\lambda$'s, then we have $\lambda_{\textmd{max}}=\textmd{max}_{i=1,\cdots,p}\lambda_i(s)\leqslant c$, where $c$ depends only on $(s_0, u_0)$.
But in the  cases of $W'$ being defined by $W'=\lambda_1\frac{\partial w_{g_1,\lambda_1}}{\partial \lambda_1}$, then $\lambda_1$  is not comparable with $\lambda_i$ for $i\ge 2$, then $\lambda_1(s)<\lambda_2\le \cdots \le \lambda_p$ since the vector field does not increase the variables $\lambda_i$ for $i\ge 2$(in fact $\lambda_i(s)=\lambda_i(0)$). In these cases, all $\lambda$'s are bounded as well. The only case is that the $\lambda$'s are comparable,  $g_i$ and $g_j$ converges to the different  critical points $\eta_i\ne\eta_j$ for all $i\ne j$. This is the case that the flow line entries the set $V_\delta(\eta_1,\cdots,\eta_p)$. If  $\eta(s, u_0)$ remains out of $V_\delta(\eta_1,\cdots,\eta_p)$,
then, $-J'(\eta(s,u_0))(\bar Z(s,u_0))\geqslant C$ since $\{\alpha_i\leqslant c, \lambda_i\leqslant c, g_i\in\mathbb{S}^{3}\}$ is a compact set. Hence, $J\left(\eta(s, u_0)\right)=J\left(\eta(0, u_0)\right)+\int_0^s J'(u)\bar Z(u)\leqslant J\left(\eta(0, u_0)\right)-C(s-s_0)$ tends to $-\infty$, when $s\rightarrow +\infty$, a
contradiction to the fact that $J$ is lower bounded.\hfill $\Box$\\
\begin{lemma}\label{ind1}
For any $u=\sum_{i=1}^p\alpha_i w_{g_i,\lambda_i}+v$ in $V_\delta(\eta_1,\cdots,\eta_p)$, we have the following expansion
of $J(u)$ after changing the variables:
$$J(u)=S\left(\sum_{i=1}^p\frac{1}{K(\eta_{\beta_i})^{\frac{2-\gamma}{\gamma}}}\right)^{\frac{\gamma}{2}}\left(1-|h|^2+\sum_{i=1}^p(|g_i^+|^2-|g_i^-|^2)+c\sum_{i=1}^p\frac{1}
{\lambda_i^2}+\|V\|^2\right),$$
where $g_i^+, g_i^-$ are the coordinates of $g_i$ near $\eta_{\beta_i}$ along the stable and unstable manifold
for K and $h=(h_1,\cdots,h_{p-1})\in \mathbb{R}^{p-1}$ with $h_i=h_i(\alpha_1,\cdots,\alpha_p)$, $i=1,\cdots, p-1$ are independent functions.
\end{lemma}
\pf
From Lemma \ref{expansion} and Theorem \ref{covering}, we have
\begin{align*}
J(u)&=\frac{\sum_{i=1}^p\alpha_i^2S}{\left(\sum_{i=1}^p\alpha_i^{\frac{4}{2-\gamma}}K(g_i)\right)^{\frac{2-\gamma}{2}}}\left[1-\frac{2-\gamma}{2}\frac{c_2}{S^{\frac{2}
{\gamma}}}\sum_{i=1}^p\frac{\alpha_i^{\frac{4}{2-\gamma}}\Delta_{\theta_1} K(g_i)}{\sum_{k=1}^p\alpha_k^{\frac{4}{2-\gamma}}K(g_k)\lambda_i^2}\right]\\&\quad
+O(\sum_{i\neq j} \varepsilon_{ij})+o(\sum_{i=1}^p\frac{1}{\lambda_i^2})+\frac{1}{2}J''
\left(\sum_{i=1}^p\alpha_iw_{
{g_i},{\lambda_i}}+\bar v(\alpha,g,\lambda)\right)V_l\cdot V_l,
\end{align*}
From the proof of Lemma \ref{Wcon}, we have $|\nabla_{\theta_1}K(g_{i})|=o(1), \varepsilon_{ij}=o(\frac{1}{\lambda_i^2})$, $|J(u)^{\frac{2-\gamma}{2}}\alpha_i^{\frac{2\gamma}{2-\gamma}}K(g_i)-1|<\varepsilon$, the expansion of the functional
$J$ can be rewritten as follows:
\begin{align*}
J(u)&=\frac{\sum_{i=1}^p\alpha_i^2}{\left(\sum_{i=1}^p\alpha_i^{\frac{4}{2-\gamma}}K(g_i)\right)^{\frac{2-\gamma}{2}}}\left[S+\frac{c}{S^{\frac{2}
{\gamma}}\sum_{k=1}^p\frac{1}{K(g_k)^{\frac{2-\gamma}{\gamma}}}}\sum_{i=1}^p\frac{-\Delta_{\theta_1} K(g_i)}{K(g_i)^{\frac{2}{\gamma}}\lambda_i^2}+o(\sum_{i=1}^p\frac{1}{\lambda_i^2})\right]\\&\quad+\frac{1}{2}J''
\left(\sum_{i=1}^p\alpha_iw_{
{g_i},{\lambda_i}}+\bar v(\alpha,g,\lambda)\right)V_l\cdot V_l.
\end{align*}
Except the term
\begin{align*}
g(\alpha,g)=\frac{\sum_{i=1}^p\alpha_i^2}{\left(\sum_{i=1}^p\alpha_i^{\frac{4}{2-\gamma}}K(g_i)\right)^{\frac{2-\gamma}{2}}},
\end{align*}
all others are positive on the right hand side of the above equality. Since $g(\alpha,g)$ is homogeneous in the variable $\alpha$, we have a degenerated critical point $(\bar \alpha_1,\cdots,\bar \alpha_p)$ which
satisfies $$\frac{\bar \alpha_i^2K(g_i)}{\bar \alpha_j^2K(g_j)}=1.$$
This critical point has an index equal to $p-1$ (since the critical point corresponds to a maximum),

On the other hand, $g(\alpha,g)$ has a single critical point $\eta = (\eta_{\beta_1}, \eta_{\beta_2}, \cdots, \eta_{\beta_p})$ in the $g$ variable. Thus, using the Morse lemma, after a change of variables, we have have the following normal form,
\begin{align*}J(u)&=\left(\sum_{i=1}^p\frac{1}{K(\eta_{\beta_i})^{\frac{2-\gamma}{\gamma}}}\right)^{\frac{\gamma}{2}}\left(S-|h|^2+\sum_{i=1}^p(|g_i^+|^2-|g_i^-|^2)+c\sum_{i=1}^p\frac{1}
{\lambda_i^2}\right)+\|V\|^2.\end{align*}\endd

For any $l$-tuple $\tau_l =(i_1, \cdots ,i_l), 1\leqslant i_j \leqslant m_1, j = 1, \cdots ,l$, let $c(\tau_l)=\left(\sum_{j=1}^l\frac{1}{K(\eta_{i_j})^{\frac{2-\gamma}{\gamma}}}\right)^{\frac{\gamma}{2}}$ denote the associated critical value. We only consider a simple situation, where for any $\tau\neq \tau', c(\tau)\neq c(\tau')$, and thus order as $c(\tau_1)<\cdots<c(\tau_{k_0})$.

From Lemma \ref{ind1} and a deformation lemma (see \cite{babook} and \cite{bahrirabino}) or directly the critical group theory (see \cite{ChangKungching93}), we have
\begin{lemma}\label{ind2}
If $c(\tau_{l-1})<a<c(\tau_l)<b<c(\tau_{l+1})$, for any coefficient group $G$, then
\begin{align*}
H_q(J_b,J_a)=
\begin{cases}
0 , \ q\neq k(\tau_l),\\
G,  \ q= k(\tau_l),
\end{cases}
\end{align*}
where $k(\tau_l)=4l-1-\sum_{j=1}^l\textmd{ind}(K,\xi_{i_j})$.
\end{lemma}
If $X$ is a topological set, then $\chi(X)$ is its Euler-Poincare characteristic with
rational coefficients.

\noindent {\it Proof of the theorem.}
Since we assumed that (\ref{main}) has no solution, $V_{\varepsilon_0}(\Sigma^+)$ is retract by deformation of $\Sigma^+$. $\Sigma^+$ is contractible, so $\chi(V_{\varepsilon_0}(\Sigma^+))=1$.
By Lemma \ref{ind2} and the Morse lemma, we have
$$\chi(V_{\varepsilon_0}(\Sigma^+))=\sum_{l=1}^{m_1}\sum_{\tau_l =(i_1, \cdots ,i_l),\eta_{i_j}\in I^+}(-1)^{4l-1-\sum_{j=1}^l\textmd{ind}(K,\eta_{i_j})}$$
is a contradiction. Therefore, (\ref{main}) has a solution $u_0\in V_{\varepsilon_0}(\Sigma^+)$.

We claim that $u_0> 0$, when $\varepsilon_0$ is small enough. Otherwise, we can write $u_0=u_0^+-u_0^-$. Multiplying equation (\ref{main}) by $u_0^-$ and integrating, using the fact that $u_0\in V_{\varepsilon_0}(\Sigma^+)$, we derive
$$||u_0^-||^2\leqslant C_1|u_0^-|_{L^{\frac{4}{2-\gamma}}}^{\frac{4}{2-\gamma}}\leqslant C_2||u_0^-||^{\frac{4}{2-\gamma}}.$$
Hence , either $u_0^- = 0$, or $||u_0^-||\geqslant C$ , where $C > 0$. Thus we have a contradiction with
$\varepsilon_0$ small enough. Therefore $u_0^- = 0$ and $u_0> 0$.

\section{Appendix}
 We first introduce some well-known inequalities which are from Taylor expansion and some computations.

\begin{lemma}
For $\alpha\geqslant 3$, there exists a constant $M>0$, such that for any $(a,b)\in \mathbb{R}^2$, there holds
\begin{align}
\left|(a+b)^{\alpha}-a^{\alpha}-\alpha a^{\alpha-1}b-\frac{\alpha(\alpha-1)}{2}a^{\alpha-2}b^2 \right|\leqslant M(|b|^{\alpha}+|a|^{\alpha-3}\textmd{inf}(|a|^3,|b|^3)).\label{ineq2}
\end{align}
\end{lemma}

In the following four lemmas, we assume $\alpha>0$.

\begin{lemma}\cite{babook}
There exists a constant $M$, such that for any $(a,b)\in \mathbb{R}^2$, $a>0,a+b>0$,
\begin{align}
\left|(a+b)^{\alpha}-a^{\alpha}-\alpha a^{\alpha-1}b\right|\leqslant M(|b|^{\alpha}+|a|^{\alpha-2}\textmd{inf}(a^2,b^2)).\label{ineq1}
\end{align}
\end{lemma}

\begin{lemma}\cite{babook}
There exists a constant $M$, such that for any $(a_1,\cdots, a_p)\in \mathbb{R}^p$,
\begin{align}
\left|\left(\sum_{i=1}^p a_i\right)^{\alpha}-\sum_{i=1}^p a_i^{\alpha}\right|\leqslant M\sum_{i\neq j}|a_i|^{\alpha-1}\textmd{inf}(|a_i|,|a_j|).\label{ineq3}
\end{align}
\end{lemma}
\begin{lemma}\cite{babook}
There exists a constant $M$, such that for any $(a_1,\cdots, a_p)\in \mathbb{R}^p$,
\begin{align}\label{ineq4}
\left|\left(\sum_{i=1}^p a_i\right)^{\alpha}-\sum_{i=1}^p a_i^{\alpha}-\alpha\sum_{i\neq j}a_i^{\alpha-1}a_j\right|\leqslant& M\left(\sum_{i\neq j}\textmd{sup}(|a_i|^{\alpha-2},|a_j|^{\alpha-2})\textmd{inf}(a_i^2,a_j^2)\right.\\&\left.\quad+\sum_{i\neq j}\textmd{inf}(|a_i|^{\alpha-1},|a_j|^{\alpha-1})\textmd{sup}(|a_i|,|a_j|)\right).\nonumber
\end{align}
\end{lemma}
\begin{lemma}\cite{babook}
There exists a constant $M$, such that for any $(a_1,\cdots, a_p)\in \mathbb{R}^p$,
\begin{align}\label{ineq5}
&\left|\left(\sum_{i=1}^p a_i\right)^{\alpha}-\sum_{i=1}^p a_i^{\alpha}-\alpha a_{i_0}^{\alpha-1}\sum_{i\neq i_0}a_i\right|\leqslant M\left(\sum_{i\neq i_0,i\neq k}|a_{i}|^{\alpha-1}\textmd{inf}(|a_i|,|a_k|)\right.\\&\left.\quad+\sum_{i\neq i_0}|a_{i_0}|^{\alpha-2}\textmd{inf}(a_{i_0}^2,a_i^2)+\sum_{i\neq i_0}\textmd{inf}(|a_{i_0}|^{\alpha-1},|a_i|^{\alpha-1})\sum_{i\neq i_0}|a_{i}|\right).\nonumber
\end{align}
\end{lemma}
\subsection{Appendix A}
We set
$$J(u)=\frac{\|u\|^2}{\left(\int_{\mathbb{S}^{3}}K(\zeta)u^{\frac{4}{2-\gamma}}\theta_1\wedge d\theta_1\right)^{\frac{2-\gamma}{2}}}=\frac{N}{D}.$$
We first expand the numerator $N$ as follows,
\begin{align*}
N&=\|u\|^2=\int_{\mathbb{S}^{3}}P_\gamma^{\theta_1} u u\theta_1\wedge d\theta_1\\&=\int_{\mathbb{S}^{3}}P_\gamma^{\theta_1} (\sum_{i=1}^p\alpha_i w_{g_i,\lambda_i}+v)(\sum_{i=1}^p\alpha_i w_{g_i,\lambda_i}+v)\theta_1\wedge d\theta_1 \\
&=\sum_{i=1}^p\alpha_i^2\int_{\mathbb{S}^{3}}P_\gamma^{\theta_1}  w_{g_i,\lambda_i} w_{g_i,\lambda_i}\theta_1\wedge d\theta_1 +\sum_{i\neq j}\alpha_i\alpha_j
\int_{\mathbb{S}^{3}}P_\gamma^{\theta_1}  w_{g_i,\lambda_i} w_{g_j,\lambda_j}\theta_1\wedge d\theta_1\\&\quad +\int_{\mathbb{S}^{3}}P_\gamma^{\theta_1} v v\theta_1\wedge d\theta_1 ,
\end{align*}
all the other terms are zero since $v$ satisfies conditions.
From now on, we denote $a_i=\mathcal{C}^{-1}(g_i),\xi=\mathcal{C}^{-1}(\zeta)$.\\

\noindent{\textbf{Lemma A.1}}
We have
$$\int_{\mathbb{S}^{3}}P_\gamma^{\theta_1}  w_{g_i,\lambda_i} w_{g_i,\lambda_i}\theta_1\wedge d\theta_1=\int_{\mathbb{H}^1}P_\gamma^{\theta_0} \delta_{a_i,\lambda_i}\delta_{a_i,\lambda_i}\theta_0\wedge d\theta_0 =S^{\frac{2}{\gamma}},$$
where $S$ is the sharp Sobolev constant given by
$$ S=\textmd{inf}_{u\in S_1^2(\mathbb{H}^1)} \frac{\int_{\mathbb{H}^1}P_\gamma u u\theta_0\wedge d\theta_0 }{(\int_{\mathbb{H}^1} u^{2^*}\theta_0\wedge d\theta_0 )^{\frac{2}{2^*}}}
,  \ 2^*=\frac{4}{2-\gamma}.$$

\noindent{\textbf{Lemma A.2}} \ It holds that for $i\neq j$
$$\int_{\mathbb{S}^{3}} w_{g_i,\lambda_i}^{\frac{2+\gamma}{2-\gamma}}w_{g_j,\lambda_j}\theta_1\wedge d\theta_1 =O'(\varepsilon_{ij}).$$
Here $O'(\hbar)$ means that when $|\hbar|< <1$, there exist two constants $C_1, C_2>0$ such that $C_1\hbar\leqslant O'(\hbar)\leqslant C_2\hbar.$

\pf  Let $\xi=(x,y,t), a_i=(x_i,y_i,t_i)\in \mathbb{H}^1$.
\begin{align*}
I&=\int_{\mathbb{S}^{3}} w_{g_i,\lambda_i}^{\frac{2+\gamma}{2-\gamma}}w_{g_j,\lambda_j}\theta_1\wedge d\theta_1\\&=\int_{\mathbb{H}^1}\delta_{a_i,\lambda_i}^{\frac{2+\gamma}{2-\gamma}}\delta_{a_j,\lambda_j}\theta_0\wedge d\theta_0 \\
&=c_0^{\frac{4}{2-\gamma}}\int_{\mathbb{H}^1}\lambda_i^{2+\gamma}\lambda_j^{2-\gamma}\frac{\theta_0\wedge d\theta_0 }{\left((1+\lambda_i^2|x-x_i|^2+\lambda_i^2|y-y_i|^2)^2+\lambda_i^4(t-t_i+2x_iy-2xy_i)^2\right)^{\frac{2+\gamma}{2}}} \\&\quad\times
\frac{1}{\left((1+\lambda_j^2|x-x_j|^2+\lambda_j^2|y-y_j|^2)^2+\lambda_j^4(t-t_j+2x_jy-2xy_j)^2\right)^{\frac{2-\gamma}{2}}}\\
&=c_0^{\frac{4}{2-\gamma}}\int_{\mathbb{H}^1}(\frac{\lambda_j}{\lambda_i})^{2-\gamma}\frac{\theta_0\wedge d\theta_0 }{\left((1+|x|^2+|y|^2)^2+t^2\right)^{\frac{2+\gamma}{2}}f_1(x,y,t)^{2-\gamma}},
\end{align*}
where
\begin{align*}
f_1(x,y,t)^2&=\left(1+(\frac{\lambda_j}{\lambda_i})^2(|x+\lambda_i(x_i-x_j)|^2+|y+\lambda_i(y_i-y_j)|^2)\right)^2\\&\quad+(\frac{\lambda_j}{\lambda_i})^4\left(t+\lambda_i^2
(t_i-t_j+2x_jy_i-2x_iy_j)+2\lambda_i(y_i-y_j)x-2\lambda_i(x_i-x_j)y\right)^2.
\end{align*}
Let $\mu$=max$(\frac{\lambda_i}{\lambda_j},\frac{\lambda_j}{\lambda_i},\lambda_i\lambda_j|d_{ij}|^2),$  where $$d_{ij}=d(a_i,a_j)=\left(((x_i-x_j)^2+(y_i-y_j)^2)^2+(t_i-t_j+2x_jy_i-2x_iy_j)^2\right)^{\frac{1}{4}}.$$
First we assume $\mu=\frac{\lambda_i}{\lambda_j}$, If $|\xi|\leqslant \frac{1}{4}\frac{\lambda_i}{\lambda_j},$ we use Taylor expansions in $\mathbb{H}^1$,
\begin{align*}
f_1(x,y,t)=g_{1}\left[1+g_{1}^{-1}Xf(0)x+g_{1}^{-1}Yf(0)y+g_{1}^{-1}T_0f(0)t+O((\frac{\lambda_j
}{\lambda_i})^2|\xi|^2)\right].
\end{align*}
Denote $$g_{1}=\left((1+\lambda_j^2|x_i-x_j|^2+\lambda_j^2|y_i-y_j|^2)^2+
\lambda_j^4(t_i-t_j+2x_jy_i-2x_iy_j)^2 \right)^{\frac{1}{2}}$$
So
\begin{align*}
f_1(x,y,t)^{-(2-\gamma)}&=g_{1}^{-(2-\gamma)}\left[1-(2-\gamma)g_{1}^{-1}Xf(0)x-(2-\gamma)g_{1}^{-1}Yf(0)y\right.\\&\quad -(2-\gamma)g_{1}^{-1}T_0f(0)\left.
t+O((\frac{\lambda_j}{\lambda_i})^2|\xi|^2)\right].
\end{align*}
Thus it yields
\begin{align*}
\int_{B(0,\frac{1}{4}\frac{\lambda_i}{\lambda_j})}\frac{\theta_0\wedge d\theta_0 }{\left((1+|x|^2+|y|^2)^2+t^2\right)^{\frac{2+\gamma}{2}}}
=\int_{\mathbb{H}^1}\frac{\theta_0\wedge d\theta_0 }{\left((1+|x|^2+|y|^2)^2+t^2\right)^{\frac{2+\gamma}{2}}}+O((\frac{\lambda_j}{\lambda_i})^{2\gamma}),
\end{align*}
\begin{align*}
\int_{B(0,\frac{1}{4}\frac{\lambda_i}{\lambda_j})}\frac{x }{\left((1+|x|^2+|y|^2)^2+t^2\right)^{\frac{2+\gamma}{2}}}\theta_0\wedge d\theta_0
=0
\end{align*}
\begin{align*}
\int_{B(0,\frac{1}{4}\frac{\lambda_i}{\lambda_j})}\frac{y }{\left((1+|x|^2+|y|^2)^2+t^2\right)^{\frac{2+\gamma}{2}}}\theta_0\wedge d\theta_0
=0
\end{align*}
\begin{align*}
\int_{B(0,\frac{1}{4}\frac{\lambda_i}{\lambda_j})}\frac{t}{\left((1+|x|^2+|y|^2)^2+t^2\right)^{\frac{2+\gamma}{2}}}\theta_0\wedge d\theta_0
=0
\end{align*}
\begin{align*}
\int_{B(0,\frac{1}{4}\frac{\lambda_i}{\lambda_j})}\frac{|\xi|^2}{\left((1+|x|^2+|y|^2)^2+t^2\right)^{\frac{2+\gamma}{2}}}\theta_0\wedge d\theta_0
=O((\frac{\lambda_i}{\lambda_j})^{2-2\gamma}),
\end{align*}
We remind that $B^c(0,\frac{1}{4}\frac{\lambda_i}{\lambda_j})=\{\xi\in \mathbb{H}^1, \mid |\xi|\geqslant \frac{1}{4}\frac{\lambda_i}{\lambda_j}\}$. Since
$f_1(x,y,t)\geqslant 1$, we have
\begin{align*}
\int_{B^c(0,\frac{1}{4}\frac{\lambda_i}{\lambda_j})}\frac{1}{\left((1+|x|^2+|y|^2)^2+t^2\right)^{\frac{2+\gamma}{2}}f(x,y,t)^{\frac{2-\gamma}{2}}}\theta_0\wedge d\theta_0
=O((\frac{\lambda_j}{\lambda_i})^{2\gamma}).
\end{align*}
Hence when $\varepsilon_{ij}$ goes to zero, we have
\begin{align*}
I&=c_0^{\frac{4}{2-\gamma}}\int_{\mathbb{H}^1}\frac{\theta_0\wedge d\theta_0 }{\left((1+|x|^2+|y|^2)^2+t^2\right)^{\frac{2+\gamma}{2}}}(\frac{\lambda_j}{\lambda_i})^{2-\gamma}g_1^{-(2-\gamma)}+O((\frac{\lambda_j}{\lambda_i})^{2+\gamma})
\\&=O'(\varepsilon_{ij}).
\end{align*}
The case  $\mu=\frac{\lambda_j}{\lambda_i}$ is similar to the case $\mu=\frac{\lambda_i}{\lambda_j}$.\\
Then we consider the third case $\mu=\lambda_i\lambda_j|d_{ij}|^2.$ In this case,
\begin{align*}
I=c_0^{\frac{4}{2-\gamma}}\int_{\mathbb{H}^1}\frac{\theta_0\wedge d\theta_0 }{\left((1+|x|^2+|y|^2)^2+t^2\right)^{\frac{2+\gamma}{2}}f_2(x,y,t)^{2-\gamma}},
\end{align*}
where
\begin{align*}
f_2(x,y,t)^2&=\left(\frac{\lambda_i}{\lambda_j}+\frac{\lambda_j}{\lambda_i}(|x+\lambda_i(x_i-x_j)|^2+|y+\lambda_i(y_i-y_j)|^2)\right)^2\\&\quad+(\frac{\lambda_j}
{\lambda_i})^2\left(t+\lambda_i^2
(t_i-t_j+2x_jy_i-2x_iy_j)+2\lambda_i(y_i-y_j)x-2\lambda_i(x_i-x_j)y\right)^2.
\end{align*}
Without loss of generality, we assume $\lambda_i\geqslant\lambda_j,$ therefore
\begin{align*}
f_2(x,y,t)=g_{2}\left[1+g_{2}^{-1}Xf(0)x+g_{2}^{-1}Yf(0)y+g_{2}^{-1}T_0f(0)t+O(\frac{\lambda_j
}{\lambda_i}g_{2}^{-1}|\xi|^2)\right].
\end{align*}
\begin{align*}
g_{2}=\left((\frac{\lambda_i}{\lambda_j}+\lambda_i\lambda_j|x_i-x_j|^2+\lambda_i\lambda_j|y_i-y_j|^2)^2+\lambda_i^2\lambda_j^2
(t_i-t_j+2x_jy_i-2x_iy_j)^2\right)^{\frac{1}{2}}.
\end{align*}
By the same arguments used in the first case, we obtain
\begin{align*}
\int_{|\xi|\leqslant\frac{\sqrt{\mu}}{10}}\frac{\theta_0\wedge d\theta_0 }{\left((1+|x|^2+|y|^2)^2+t^2\right)^{\frac{2+\gamma}{2}}f_2(x,y,t)^{2-\gamma}}
=O'(\varepsilon_{ij}).
\end{align*}
Let
\begin{align*}
B_1:=\{\xi\in \mathbb{H}^1 \mid |\xi+\lambda_id_{ij}|\leqslant\frac{1}{10}\lambda_id_{ij}\},
\end{align*}
\begin{align*}
B_2:=\{\xi\in \mathbb{H}^1 \mid |\xi|\leqslant\frac{\sqrt{\mu}}{10}\}.
\end{align*}
We have
\begin{align*}
\int_{(B_1\cup B_2)^c} \frac{\theta_0\wedge d\theta_0 }{\left((1+|x|^2+|y|^2)^2+t^2\right)^{\frac{2+\gamma}{2}}f_2(x,y,t)^{2-\gamma}}=O(\varepsilon_{ij}^{\frac{2}{2-\gamma}}).
\end{align*}
On $B_1$, we have $|\xi|\geqslant\frac{9}{10}\lambda_i|d_{ij}|$ , we obtain
\begin{align*}
\int_{B_1} \frac{\theta_0\wedge d\theta_0 }{\left((1+|x|^2+|y|^2)^2+t^2\right)^{\frac{2+\gamma}{2}}f_2(x,y,t)^{2-\gamma}}=O(\varepsilon_{ij}^{\frac{2}{2-\gamma}}).
\end{align*}
This completes the proof. \endd

Similar to the proof of Lemma A.2, we have the following results.

\noindent{\textbf{Lemma A.3}} (1)\ It holds that
$$\int_{\mathbb{S}^{3}} w_{g_i,\lambda_i}^{\frac{2}{2-\gamma}}w_{g_j,\lambda_j}^{\frac{2}{2-\gamma}}\theta_1\wedge d\theta_1=\int_{\mathbb{H}^1}\delta_{a_i,\lambda_i}^{\frac{2}{2-\gamma}}\delta_{a_j,\lambda_j}^{\frac{2}{2-\gamma}}\theta_0\wedge d\theta_0 =O(\varepsilon_{ij}^
{\frac{2}{2-\gamma}}\ln\varepsilon_{ij}^{-1}).$$\\
(2)\ Let $\alpha, \beta > 1$, such that $\alpha+\beta=\frac{4}{2-\gamma}, \theta=$ inf $(\alpha, \beta)$, it holds that
$$\int_{\mathbb{S}^{3}} w_{g_i,\lambda_i}^\alpha w_{g_j,\lambda_j}^\beta\theta_1\wedge d\theta_1=\int_{\mathbb{H}^1}\delta_{a_i,\lambda_i}^\alpha\delta_{a_j,\lambda_j}^\beta\theta_0\wedge d\theta_0 =O(\varepsilon_{ij}^
\theta(\ln\varepsilon_{ij}^{-1})^{\frac{2-\gamma}{2}\theta}).$$\\

Let us consider the denominator $D$ of $J$,
$$D^{\frac{2}{2-\gamma}}=\int_{\mathbb{S}^{3}}K(\zeta)\left(\sum_{i=1}^p\alpha_i w_{g_i,\lambda_i}+v\right)^{\frac{4}{2-\gamma}}\theta_1\wedge d\theta_1 .$$

\noindent{\textbf{Lemma A.4}} If $|v|<\sum_{i=1}^p\alpha_i w_{g_i,\lambda_i},$
\begin{align}\lb{A.4}
\nonumber\int_{\mathbb{S}^3}K(\zeta)\left(\sum_{i=1}^p\alpha_i w_{g_i,\lambda_i}+v\right)^{\frac{4}{2-\gamma}}&=
\int_{\mathbb{S}^3}K(\zeta)\left(\sum_{i=1}^p\alpha_i w_{g_i,\lambda_i}\right)^{\frac{4}{2-\gamma}}\\\nonumber &\quad+
\frac{4}{2-\gamma}\int_{\mathbb{S}^3}K(\zeta)\left(\sum_{i=1}^p\alpha_iw_{g_i,\lambda_i}\right)^{\frac{2+\gamma}{2-\gamma}}v\\\nonumber &\quad+
\frac{2(2+\gamma)}{(2-\gamma)^2}\int_{\mathbb{S}^3}K(\zeta)\left(\sum_{i=1}^p\alpha_i w_{g_i,\lambda_i}\right)^{\frac{2\gamma}{2-\gamma}}v^2\\&\quad+
O(\|v\|^3).
\end{align}
\pf Using (\ref{ineq2}), we have
\begin{align*}
&\left|\int_{\mathbb{S}^3}K(\zeta)\left(\sum_{i=1}^p\alpha_i w_{g_i,\lambda_i}+v\right)^{\frac{4}{2-\gamma}}-\right.
\int_{\mathbb{S}^3}K(\zeta)\left(\sum_{i=1}^p\alpha_i w_{g_i,\lambda_i}\right)^{\frac{4}{2-\gamma}}\\&\quad\left.-
\frac{4}{2-\gamma}\int_{\mathbb{S}^3}K(\zeta)\left(\sum_{i=1}^p\alpha_i w_{g_i,\lambda_i}\right)^{\frac{2+\gamma}{2-\gamma}}v-
\frac{2(2+\gamma)}{(2-\gamma)^2}\int_{\mathbb{S}^3}K(\zeta)\left(\sum_{i=1}^p\alpha_i w_{g_i,\lambda_i}\right)^{\frac{2\gamma}{2-\gamma}}v^2\right| \\&\leqslant
M_2\left(\int_{\mathbb{S}^3}|v|^{\frac{4}{2-\gamma}}+\int_{\mathbb{S}^3}\left(\sum_{i=1}^p\alpha_i w_{g_i,\lambda_i}\right)^{\frac{3\gamma-2}{2-\gamma}}\textmd
{inf}\left((\sum_{i=1}^p
\alpha_i w_{g_i,\lambda_i})^3,|v|^3\right)\right)\\&=O(\|v\|^3).
\end{align*}
\endd

In order to get more information from Lemma A.4, we now estimate the first three terms on the right hand side of \eqref{A.4}.

\noindent{\textbf{Lemma A.4.1}}
\begin{align*}
&\int_{\mathbb{S}^3}K(\zeta)(\sum_{i=1}^p\alpha_i w_{g_i,\lambda_i})^{\frac{4}{2-\gamma}}\theta_1\wedge d\theta_1\\ &=
\sum_{i=1}^p\alpha_i^{\frac{4}{2-\gamma}}\left(K(g_i)S^{\frac{2}{\gamma}}+\frac{c_2\Delta_{\theta_1} K(g_i)}{\lambda_i^2}\right)
+O'(\sum_{i\neq j}\varepsilon_{ij})+o(\frac{1}{\lambda_i^2}),
\end{align*}
where $c_2=4c_0^{\frac{4}{2-\gamma}}\int_{\mathbb{H}^1}\frac{(x^2+y^2)\theta_0\wedge d\theta_0 }{\left((1+|x|^2+|y|^2)^2+t^2\right)^2}.$\\
\pf From definition and \eqref{relat2}, there holds
 \begin{align*}
\int_{\mathbb{S}^3}K(\zeta)(\sum_{i=1}^p\alpha_i w_{g_i,\lambda_i})^{\frac{4}{2-\gamma}}\theta_1\wedge d\theta_1=\int_{\mathbb{H}^1}\tilde{K}(\xi)\left(\sum_{i=1}^p\alpha_i\delta_{a_i,\lambda_i}\right)^{\frac{4}{2-\gamma}}\theta_0\wedge d\theta_0,
\end{align*}
where $\tilde{K}=K\circ \mathcal{C}$.
From (\ref{ineq4}) and Lemma A.3, we get
\begin{align*}
&\left| \int_{\mathbb{H}^1}\tilde{K}(\xi)\left(\sum_{i=1}^p\alpha_i\delta_{a_i,\lambda_i}\right)^{\frac{4}{2-\gamma}}\theta_0\wedge d\theta_0\right.\\&\quad -
\left.\int_{\mathbb{H}^1}\tilde{K}(\xi)\left(\sum_{i=1}^p(\alpha_i\delta_{a_i,\lambda_i})^{\frac{4}{2-\gamma}}+\frac{4}{2-\gamma}\sum_{i\neq j}(\alpha_i\delta_{a_i,\lambda_i})^{\frac{2+\gamma}{2-\gamma}}\alpha_j\delta_{a_j,\lambda_j}\right)\theta_0\wedge d\theta_0\right|\\&\leqslant M_4 \sum_{i\neq j}\int_{\mathbb{H}^1}
\textmd{sup}\left((\alpha_i\delta_{a_i,\lambda_i})^{\frac{2\gamma}{2-\gamma}},(\alpha_j\delta_{a_j,\lambda_j})^{\frac{2\gamma}{2-\gamma}}\right)\textmd{inf}\left((\alpha_i\delta_{a_i,
\lambda_i})^2,(\alpha_j\delta_{a_j,\lambda_j})^2\right)
\theta_0\wedge d\theta_0 \\ &\quad+M_4 \sum_{i\neq j}\int_{\mathbb{H}^1}
\textmd{inf}\left((\alpha_i\delta_{a_i,\lambda_i})^{\frac{2+\gamma}{2-\gamma}},(\alpha_j\delta_{a_j,\lambda_j})^{\frac{2+\gamma}{2-\gamma}}\right)\textmd{sup}\left(\alpha_i
\delta_{a_i,\lambda_i},\alpha_j\delta_{a_j,\lambda_j}\right)
\theta_0\wedge d\theta_0\\&
\leqslant 2 M_4\sum_{i\neq j}(\alpha_i\alpha_j)^{\frac{2}{2-\gamma}}\int_{\mathbb{H}^1}\delta_{a_i,\lambda_i}^{\frac{2}{2-\gamma}}\delta_{a_j,\lambda_j}^{\frac{2}{2-\gamma}}
\theta_0\wedge d\theta_0 \\&=O(\sum_{i\neq j}\varepsilon_{ij}^
{\frac{2}{2-\gamma}}\ln\varepsilon_{ij}^{-1}).
\end{align*}

Then we have
\begin{align*}
\int_{\mathbb{H}^1}\tilde{K}(\xi)\delta_{a_i,\lambda_i}^{\frac{4}{2-\gamma}}\theta_0\wedge d\theta_0 &=
\int_{\mathbb{H}^1}\left(\tilde{K}(\xi)-\tilde{K}(a_i)\right)\delta_{a_i,\lambda_i}^{\frac{4}{2-\gamma}}\theta_0\wedge d\theta_0 +\tilde{K}(a_i)S^{\frac{2}{\gamma}}\\
&=\tilde{K}(a_i)S^{\frac{2}{\gamma}}+\int_{B(a_i,\varepsilon)}\left(\tilde{K}(\xi)-\tilde{K}(a_i)\right)\delta_{a_i,\lambda_i}^{\frac{4}{2-\gamma}}\\&\quad+\int_{B^c(a_i,\varepsilon)}\left(
\tilde{K}(\xi)-\tilde{K}(a_i)\right)
\delta_{a_i,\lambda_i}^{\frac{4}{2-\gamma}}
\end{align*}
and
\begin{align*}
\int_{B^c(a_i,\varepsilon)}\delta_{a_i,\lambda_i}^{\frac{4}{2-\gamma}}\theta_0\wedge d\theta_0 =O(\frac{1}{\lambda_i^{4}}).
\end{align*}
By Taylor expansion, there holds
\begin{align*}
\int_{B(a_i,\varepsilon)}(\tilde{K}(\xi)-\tilde{K}(a_i))\delta_{a_i,\lambda_i}^{\frac{4}{2-\gamma}}\theta_0\wedge d\theta_0 =\frac{c_2\Delta_{\theta_0} \tilde{K}(a_i)}{4\lambda_i^2}+o(\frac{1}{\lambda_i^2}).
\end{align*}
Finally we estimate for $i\neq j$, by Lemma A.2, Lemma A.3 and Taylor expansion and Young inequality,
\begin{align*}
&\int_{\mathbb{H}^1}\tilde{K}(\xi)(\delta_{a_i,\lambda_i})^{\frac{2+\gamma}{2-\gamma}}\delta_{a_j,\lambda_j}\theta_0\wedge d\theta_0
\\&=\int_{\mathbb{H}^1}\tilde{K}(a_i)(\delta_{a_i,\lambda_i})^{\frac{2+\gamma}{2-\gamma}}\delta_{a_j,\lambda_j}+
\int_{B(a_i,\varepsilon)}(\tilde{K}(\xi)-\tilde{K}(a_i))(\delta_{a_i,\lambda_i})^{\frac{2+\gamma}{2-\gamma}}\delta_{a_j,\lambda_j}
\\&\quad+\int_{B^c(a_i,\varepsilon)}(\tilde{K}(\xi)-\tilde{K}(a_i))(\delta_{a_i,\lambda_i})^{\frac{2+\gamma}{2-\gamma}}\delta_{a_j,\lambda_j}\\&
=O(\varepsilon_{ij})+
O'(\varepsilon_{ij}^{\frac{2}{2-\gamma}}\ln\varepsilon_{ij}^{-1})+o(\frac{1}{\lambda_i^2})+O(\frac{1}{\lambda_i^{2+\gamma}\lambda_j^{2-\gamma}}).
\end{align*}\endd

\noindent{\textbf{Lemma A.4.2}} We have
\begin{align*}
\int_{\mathbb{S}^3}K(\zeta)(\sum_{i=1}^p\alpha_i w_{g_i,\lambda_i})^{\frac{2+\gamma}{2-\gamma}}v=
\sum_{i=1}^p O\left(\frac{|\nabla_{\theta_1} K(g_i)|}{\lambda_i}+\frac{1}{\lambda_i^2}\right)\|v\|+O(\sum_{i\neq j}\varepsilon_{ij}(\ln\varepsilon_{ij}^{-1})^{{\frac{2-\gamma}{2}}})\|v\|.
\end{align*}
\pf By (\ref{relat2}), we have
\begin{align*}
\int_{\mathbb{S}^3}K(\zeta)(\sum_{i=1}^p\alpha_i w_{g_i,\lambda_i})^{\frac{2+\gamma}{2-\gamma}}v\theta_1\wedge d\theta_1\leqslant C_0\int_{\mathbb{H}^1}\tilde{K}(\xi)\left(\sum_{i=1}^p\alpha_i\delta_{a_i,\lambda_i}\right)^{\frac{2+\gamma}{2-\gamma}}v\theta_0\wedge d\theta_0,
\end{align*}
where $C_0$ is a constant depending on $\gamma$.
From (\ref{ineq3}) , H\"{o}lder inequality and Lemma A.3, we get
\begin{align*}
&\left|\int_{\mathbb{H}^1}\tilde{K}(\xi)\left(\sum_{i=1}^p\alpha_i\delta_{a_i,\lambda_i}\right)^{\frac{2+\gamma}{2-\gamma}}v\theta_0\wedge d\theta_0-
\sum_{i=1}^p\alpha_i^{\frac{2+\gamma}{2-\gamma}}\int_{\mathbb{H}^1}\tilde{K}(\xi)(\delta_{a_i,\lambda_i})^{\frac{2+\gamma}{2-\gamma}}v\right|
\\&\leqslant M_3 \sum_{i\neq j}\int_{\mathbb{H}^1}(\alpha_i\delta_{a_i,\lambda_i})^{\frac{2\gamma}{2-\gamma}}\textmd{inf}(\alpha_i\delta_{a_i,\lambda_i},
\alpha_j\delta_{a_j,\lambda_j})|v|
\\&\leqslant M_3 \sum_{i\neq j}\left(\int_{\mathbb{H}^1}[(\alpha_i\delta_{a_i,\lambda_i})^{\frac{2\gamma}{2-\gamma}}
\alpha_j\delta_{a_j,\lambda_j}]^{\frac{4}{2+\gamma}}\right)^\frac{2+\gamma}{4}\|v\|
\\&=M_3 \sum_{i\neq j}\left(\int_{\mathbb{H}^1}[(\alpha_i\delta_{a_i,\lambda_i})^{\frac{4(3\gamma-2)}{(2-\gamma)(2+\gamma)}}
(\alpha_i\delta_{a_i,\lambda_i}\alpha_j\delta_{a_j,\lambda_j})^{\frac{2+\gamma}{4}}]^{\frac{4}{2+\gamma}}\right)^\frac{2+\gamma}{4}\|v\|
\\&\leqslant M_3 \sum_{i\neq j}\left(\left(\int_{\mathbb{H}^1}(\alpha_i\delta_{a_i,\lambda_i})^{\frac{4}{2-\gamma}}\right)^{\frac{3\gamma-2}{2+\gamma}}\left(\int_{\mathbb{H}^1}
(\alpha_i\delta_{a_i,\lambda_i}\alpha_j\delta_{a_j,\lambda_j})^{\frac{4}{2-\gamma}}\right)^{\frac{4-2\gamma}{2+\gamma}}\right)^\frac{2+\gamma}{4}\|v\|
\\&=O(\varepsilon_{ij}(\ln\varepsilon_{ij}^{-1})^{{\frac{2-\gamma}{2}}})\|v\|.
\end{align*}

Since $v$ satisfies (\ref{V0}), we have
\begin{align*}
\int_{\mathbb{H}^1}\tilde{K}(\xi)(\delta_{a_i,\lambda_i})^{\frac{2+\gamma}{2-\gamma}}v&=\int_{\mathbb{H}^1}(\tilde{K}(\xi)-\tilde{K}(a_i))(\delta_{a_i,\lambda_i})^{\frac{2+\gamma}
{2-\gamma}}v\\&=\int_{B(a_i,\varepsilon)}\left(\nabla_{\theta_0} \tilde{K}(a_i)(\xi-a_i)+O(|\xi-a_i|^2)\right)(\delta_{a_i,\lambda_i})^{\frac{2+\gamma}{2-\gamma}}v
\\&\quad+\int_{B^c(a_i,\varepsilon)}(\tilde{K}(\xi)-\tilde{K}(a_i))(\delta_{a_i,\lambda_i})^{\frac{2+\gamma}{2-\gamma}}v\\&\leqslant
C\left(\frac{|\nabla_{\theta_0} \tilde{K}(a_i)|}{\lambda_i}+\frac{1}{\lambda_i^2}+\frac{1}{\lambda_i^{2+\gamma}}\right)\|v\|.
\end{align*} \endd

\noindent{\textbf{Lemma A.4.3}} We have
\begin{align*}
&\int_{\mathbb{S}^3}K(\zeta)(\sum_{i=1}^p\alpha_i w_{g_i,\lambda_i})^{\frac{2\gamma}{2-\gamma}}v^2 =
\sum_{i=1}^p \alpha_i^{\frac{2\gamma}{2-\gamma}}K(g_i)\int_{\mathbb{S}^3} w_{g_i,\lambda_i}^{\frac{2\gamma}{2-\gamma}}v^2\\&\quad\quad\quad\quad\quad\quad\quad\quad\quad\quad\quad\quad+O(\|v\|^2)\left(\sum_{i=1}^p\frac{|\nabla_{\theta_1} K(g_i)|}{\lambda_i}+\sum_{i=1}^p\frac{1}{\lambda_i^{2\gamma}}+\sum_{i\neq j}
\varepsilon_{ij}^{\frac{\gamma}{2-\gamma}}(\ln\varepsilon_{ij}^{-1})^{\frac{\gamma}{2}}\right).
\end{align*}
\pf
Using (\ref{ineq3}), we get
\begin{align*}
&\left|\int_{\mathbb{S}^3}K(\zeta)(\sum_{i=1}^p\alpha_i w_{g_i,\lambda_i})^{\frac{2\gamma}{2-\gamma}}v^2\theta_1\wedge d\theta_1 -
\int_{\mathbb{S}^3}K(\zeta)\sum_{i=1}^p(\alpha_i w_{g_i,\lambda_i})^{\frac{2\gamma}{2-\gamma}}v^2\right| \\&\quad\leqslant M_3\sum_{i\neq j}\int_{\mathbb{S}^3}(\alpha_i\delta_{a_i,\lambda_i})^{\frac{3\gamma-2}{2-\gamma}}\textmd{inf}(\alpha_i\delta_{a_i,\lambda_i},
\alpha_j\delta_{a_j,\lambda_j})v^2.
\end{align*}
By H\"{o}lder inequality and Lemma A.3, we can easily get,
$$\sum_{i\neq j}\int_{\mathbb{S}^3}(\alpha_i\delta_{a_i,\lambda_i})^{\frac{3\gamma-2}{2-\gamma}}\textmd{inf}(\alpha_i\delta_{a_i,\lambda_i},
\alpha_j\delta_{a_j,\lambda_j})v^2 =O\left(\sum_{i\neq j}
\varepsilon_{ij}^{\frac{\gamma}{2-\gamma}}(\ln\varepsilon_{ij}^{-1})^{\frac{\gamma}{2}}\|v\|^2\right).$$
Now, we compute
\begin{align*}
&\int_{\mathbb{S}^3}K(\zeta)(\alpha_i w_{g_i,\lambda_i})^{\frac{2\gamma}{2-\gamma}}v^2\theta_1\wedge d\theta_1\\&=\int_{\mathbb{S}^3}K(g_i)(\alpha_i w_{g_i,\lambda_i})^{\frac{2\gamma}{2-\gamma}}v^2\theta_1\wedge d\theta_1+\int_{\mathbb{S}^3}(K(\zeta)-K(g_i))(\alpha_i w_{g_i,\lambda_i})^{\frac{2\gamma}{2-\gamma}}v^2\theta_1\wedge d\theta_1\\ &\leqslant
\int_{\mathbb{S}^3}K(g_i)(\alpha_i w_{g_i,\lambda_i})^{\frac{2\gamma}{2-\gamma}}v^2\theta_1\wedge d\theta_1+\int_{\mathbb{H}^1}(\tilde{K}(\xi)-\tilde{K}(a_i))(\delta_{a_i,\lambda_i})^{\frac{2\gamma}{2-\gamma}}
v^2\\&= \int_{\mathbb{S}^3}K(g_i)(\alpha_i w_{g_i,\lambda_i})^{\frac{2\gamma}{2-\gamma}}v^2\theta_1\wedge d\theta_1+\int_{B^c(a_i,\varepsilon)}(\tilde{K}(\xi)-\tilde{K}(a_i))(\delta_{a_i,\lambda_i})^{\frac{2\gamma}{2-\gamma}}v^2\\&\quad+\int_{B(a_i,\varepsilon)}(\tilde{K}(\xi)-\tilde{K}(a_i))(\delta_{a_i,\lambda_i})^{\frac{2\gamma}{2-\gamma}}v^2
\\&=\int_{\mathbb{S}^3}K(g_i)(\alpha_i w_{g_i,\lambda_i})^{\frac{2\gamma}{2-\gamma}}v^2\theta_1\wedge d\theta_1+O\left(\frac{1}{\lambda_i^{2\gamma}}\right)\|v\|^2\\&\quad+\int_{B(a_i,\varepsilon)}\left(\nabla_{\theta_0} \tilde{K}(a_i)(\xi-a_i)+O(|\xi-a_i|^2)\right)(\delta_{a_i,\lambda_i})^{\frac{2\gamma}{2-\gamma}}v^2
\\&=\int_{\mathbb{S}^3}K(g_i)(\alpha_i w_{g_i,\lambda_i})^{\frac{2\gamma}{2-\gamma}}v^2\theta_1\wedge d\theta_1+O\left(\frac{1}{\lambda_i^{2\gamma}}\right)\|v\|^2+O\left(\frac{|\nabla_{\theta_1} K(g_i)|}{\lambda_i}\right)\|v\|^2.
\end{align*}\endd

Now we can complete the Proof of Lemma \ref{expansion}.
\begin{align*}
D^{\frac{2}{2-\gamma}}&=\sum_{i=1}^p\alpha_i^{\frac{4}{2-\gamma}}\left(K(g_i)S^{\frac{2}{\gamma}}+\frac{c_2\Delta_{\theta_1} K(g_i)}{\lambda_i^2}\right)
+\frac{4}{2-\gamma}\int_{\mathbb{S}^3}K(\zeta)\left(\sum_{i=1}^p\alpha_i w_{g_i,\lambda_i}\right)^{\frac{2+\gamma}{2-\gamma}}v\\&\quad+
\frac{2(2+\gamma)}{(2-\gamma)^2}\int_{\mathbb{S}^3}K(\zeta)\left(\sum_{i=1}^p\alpha_i w_{g_i,\lambda_i}\right)^{\frac{2\gamma}{2-\gamma}}v^2+O'(\sum_{i\neq j}\varepsilon_{ij})+o(\frac{1}{\lambda_i^2})+O(\|v\|^3)
\\&=(\sum_{i=1}^p\alpha_i^{\frac{4}{2-\gamma}}K(g_i)S^{\frac{2}{\gamma}})\left[1+\left(\sum_{i=1}^p\alpha_i^{\frac{4}{2-\gamma}}K(g_i)S^{\frac{2}{\gamma}}\right)^{-1}\frac{c_2
\Delta_{\theta_1} K(g_i)}{\lambda_i^2}\right.\\&\quad+(\sum_{i=1}^p\alpha_i^{\frac{4}{2-\gamma}}K(g_i)S^{\frac{2}{\gamma}})^{-1}\frac{4}{2-\gamma}\int_{\mathbb{S}^3}K(\zeta)\left(\sum_{i=1}^p\alpha_i w_{g_i,\lambda_i}\right)^{\frac{2+\gamma}{2-\gamma}}v\\&\quad+(\sum_{i=1}^p\alpha_i^{\frac{4}{2-\gamma}}K(g_i)S^{\frac{2}{\gamma}})^{-1}\int_{\mathbb{S}^3}K(\zeta)\left(\sum_{i=1}^p\alpha_i w_{g_i,\lambda_i}\right)^{\frac{2\gamma}{2-\gamma}}v^2 \\&\left.\quad+O'(\sum_{i\neq j}\varepsilon_{ij})+o(\frac{1}{\lambda_i^2})+O(\|v\|^3)\right].
\end{align*}
Therefore combining Lemma A.1, Lemma A.2 and Lemma A.3, we have
\begin{align*}
J(u)&=\frac{\sum_{i=1}^p\alpha_i^2S}{\left(\sum_{i=1}^p\alpha_i^{\frac{4}{2-\gamma}}K(g_i)\right)^{\frac{2-\gamma}{2}}}\left[1-\frac{2-\gamma}{2}\frac{c_2}{S^{\frac{2}
{\gamma}}}\sum_{i=1}^p\frac{\alpha_i^{\frac{4}{2-\gamma}}\Delta_{\theta_1} K(g_i)}{\sum_{k=1}^p\alpha_k^{\frac{4}{2-\gamma}}K(g_k)\lambda_i^2}\right.\\&\quad
-\frac{2}{\sum_{k=1}^p\alpha_k^{\frac{4}{2-\gamma}}K(g_k)S^{\frac{2}{\gamma}}}\int_{\mathbb{S}^3}K(\zeta)\left(\sum_{i=1}^p\alpha_iw_{g_i,\lambda_i}\right)^
{\frac{2+\gamma}{2-\gamma}}v\\&\quad+\frac{1}{\sum_{k=1}^p\alpha_k^2S^{\frac{2}{\gamma}}}\|v\|^2-\frac{2+\gamma}{(2-\gamma)S^{\frac{2}{\gamma}}\sum_{k=1}^p\alpha_k^{\frac{4}
{2-\gamma}}K(g_k)}\int_{\mathbb{S}^3}K(\zeta)\sum_{i=1}^p(\alpha_iw_{g_i,\lambda_i})^{\frac{2\gamma}{2-\gamma}}v^2\\&\quad\left.+O'(\sum_{i\neq j} \varepsilon_{ij})+o(\sum_{i=1}^p\frac{1}{\lambda_i^2})+o(\|v\|^2)\right].
\end{align*}
This is the estimate of $J(u)$ in Lemma \ref{expansion}.

\subsection{Appendix B}

We define
$$\lambda(u)=\left(\int_{\mathbb{S}^3}K(\zeta)u^{\frac{4}{2-\gamma}}\theta_1\wedge d\theta_1 \right)^{\frac{-(2-\gamma)}{2}},$$ so
$$J(u)=\lambda(u)\int_{\mathbb{S}^3}P_\gamma uu\,\theta_1\wedge d\theta_1.$$
By direct computation,
$$\lambda'(u)W=-2\lambda(u)^{\frac{4-\gamma}{2-\gamma}}\int_{\mathbb{S}^3}K(\zeta)u^{\frac{2+\gamma}{2-\gamma}}W\theta_1\wedge d\theta_1 .$$
Thus we have
\begin{align*}
J'(u)W&=\lambda'(u)W\int_{\mathbb{S}^3}P_\gamma uu\theta_1\wedge d\theta_1+2\lambda(u)\int_{\mathbb{S}^3}P_\gamma uW\theta_1\wedge d\theta_1\\&=
2\lambda(u)\left(\int_{\mathbb{S}^3}P_\gamma uW\theta_1\wedge d\theta_1-\lambda(u)^{\frac{2}{2-\gamma}}\int_{\mathbb{S}^3}K(\zeta)u^{\frac{2+\gamma}{2-\gamma}}W\theta_1\wedge d\theta_1\int_{\mathbb{S}^3}P_\gamma uu\theta_1\wedge d\theta_1\right).
\end{align*}
Since for $u_0=\sum_{i=1}^p\alpha_i w_{g_i,\lambda_i}\in V(p,\varepsilon)\subset \Sigma^+$, we have $\int_{\mathbb{S}^3}P_\gamma u_0u_0\,\theta_1\wedge d\theta_1 =1$.
Thus
\begin{align}\label{J'}
J'(u_0)W=
2\lambda(u_0)\left(\int_{\mathbb{S}^3}P_\gamma u_0W-\lambda(u_0)^{\frac{2}{2-\gamma}}\int_{\mathbb{S}^3}K(\zeta)u_0^{\frac{2+\gamma}{2-\gamma}}W\right).
\end{align}

We first take $W=\lambda_j\frac{\partial w_{g_j,\lambda_j}}{\partial\lambda_j}$ in (\ref{J'}), we obtain
\begin{align*}
J'(u)\left(\lambda_j\frac{\partial w_{g_j,\lambda_j}}{\partial\lambda_j}\right)&=
2\lambda(u)\left(\int_{\mathbb{S}^3}P_\gamma \left(\sum_{i=1}^p\alpha_i w_{g_i,\lambda_i}\right) \left(\lambda_j\frac{\partial w_{g_j,\lambda_j}}{\partial\lambda_j}\right)\right.\\&\left.\quad-\lambda(u)^{\frac{2}{2-\gamma}}\int_{\mathbb{S}^3}K(\zeta)\left(\sum_{i=1}^p\alpha_i w_{g_i,
\lambda_i}\right)^
{\frac{2+\gamma}{2-\gamma}}\left(\lambda_j\frac{\partial w_{g_j,\lambda_j}}{\partial\lambda_j}\right)\right).
\end{align*}

In the remainder of the part B, we will give some lemmas to complete the proof of  (\ref{es01}).

\noindent{\textbf{Lemma B.1}} We have
\begin{align*}
\int_{\mathbb{S}^3}P_\gamma ( w_{g_i,\lambda_i})\lambda_i\frac{\partial w_{g_i,\lambda_i}}{\partial\lambda_i}\theta_1\wedge d\theta_1 =0.
\end{align*}
\pf Since $\int_{\mathbb{S}^3}P_\gamma ( w_{g_i,\lambda_i}) w_{g_i,\lambda_i}\theta_1\wedge d\theta_1=S^{\frac{2}{\gamma}}$ is independent of $\lambda_i$, we get the result.\endd

\noindent{\textbf{Lemma B.2}}
For $i\neq j$, we have
\begin{align*}
\int_{\mathbb{S}^3}P_\gamma ( w_{g_i,\lambda_i})\lambda_j\frac{\partial w_{g_j,\lambda_j}}{\partial\lambda_j}\theta_1\wedge d\theta_1
=O'(\lambda_j\frac{\partial\varepsilon_{ij}}{\partial\lambda_j}).
\end{align*}
In the proof of this result,  the idea is same as that in the proof of Lemma A.2. The details are omitted.

\noindent{\textbf{Lemma B.3}} There holds
\begin{align*}
&\int_{\mathbb{S}^3}K(\zeta)\left(\sum_{i=1}^p\alpha_i w_{g_i,\lambda_i}\right)^
{\frac{2+\gamma}{2-\gamma}}\left(\lambda_j\frac{\partial w_{g_j,\lambda_j}}{\partial\lambda_j}\right)\\&=-\alpha_j^
{\frac{2+\gamma}{2-\gamma}}\frac{2-\gamma}{4}c_2\frac{\Delta_{\theta_1} K(g_j)}{\lambda_j^2}+\sum_{i\neq j}O'(\lambda_j\frac{\partial\varepsilon_{ij}}{\partial\lambda_j})+
o(\sum_{i\neq j} \varepsilon_{ij})+o\left(\frac{1}{\lambda_j^2}\right).
\end{align*}
\pf Using (\ref{ineq5}), we have
\begin{align*}
&\left|\int_{\mathbb{S}^3}K(\zeta)\left(\sum_{i=1}^p\alpha_i w_{g_i,\lambda_i}\right)^
{\frac{2+\gamma}{2-\gamma}}\left(\lambda_j\frac{\partial w_{g_j,\lambda_j}}{\partial\lambda_j}\right)\right.\\&\quad\left.-\int_{\mathbb{S}^3}K(\zeta)\left[
\sum_{i=1}^p(\alpha_i w_{g_i,\lambda_i})^
{\frac{2+\gamma}{2-\gamma}}+\frac{2+\gamma}{2-\gamma}(\alpha_j w_{g_j,\lambda_j})^
{\frac{2\gamma}{2-\gamma}}\left(\sum_{i\neq j}\alpha_i w_{g_i,\lambda_i}\right)\right]\lambda_j\frac{\partial w_{g_j,\lambda_j}}{\partial\lambda_j}\right|\\&\leqslant M_5\left(\int_{\mathbb{S}^3}K(\zeta)\sum_{i\neq j,i\neq k}(\alpha_i w_{g_i,\lambda_i})^{\frac{2\gamma}{2-\gamma}}\textmd{inf}(\alpha_i w_{g_i,\lambda_i},
\alpha_k w_{g_k,\lambda_k})\lambda_j\frac{\partial w_{g_j,\lambda_j}}{\partial\lambda_j}\right)\\&\quad+ M_5\left(\int_{\mathbb{S}^3}K(\zeta)\sum_{i\neq j}(\alpha_j w_{g_j,\lambda_j})^
{\frac{3\gamma-2}{2-\gamma}}\textmd{inf}((\alpha_i w_{g_i,\lambda_i})^2,
(\alpha_j w_{g_j,\lambda_j})^2)\lambda_j\frac{\partial w_{g_j,\lambda_j}}{\partial\lambda_j}\right).
\end{align*}
Together with Lemma B.3.1, Lemma B.3.2, Lemma B.3.3, Lemma B.3.4 and Lemma B.3.5 below, we get a complete proof.\endd

\noindent{\textbf{Lemma B.3.1}} We have
\begin{align*}
\int_{\mathbb{S}^3}K(\zeta) w_{g_i,\lambda_i}^
{\frac{2+\gamma}{2-\gamma}}\lambda_j\frac{\partial w_{g_j,\lambda_j}}{\partial\lambda_j}\theta_1
\wedge d\theta_1=-\frac{2-\gamma}{2}c_2\frac{\Delta_{\theta_1} K(g_i)}{\lambda_i^2}+o\left(\frac{1}{\lambda_i^2}\right).
\end{align*}
\pf There hold
\begin{align*}
&\int_{\mathbb{S}^3}K(\zeta) w_{g_i,\lambda_i}^
{\frac{2+\gamma}{2-\gamma}}\lambda_j\frac{\partial w_{g_j,\lambda_j}}{\partial\lambda_j}\theta_1
\wedge d\theta_1\\&=\int_{\mathbb{H}^1}\tilde{K}(\xi)(\delta_{a_i,\lambda_i})^{\frac{2+\gamma}{2-\gamma}}\lambda_i\frac{\partial\delta_{a_i,\lambda_i}}{\partial\lambda_i}\theta_0\wedge d
\theta_0 \\&=\frac{2-\gamma}{4}\lambda_i\frac{\partial}{\partial\lambda_i}\int_{\mathbb{H}^1}\tilde{K}(\xi)(\delta_{a_i,\lambda_i})^{\frac{4}{2-\gamma}}\theta_0
\wedge d\theta_0 \\&=
\frac{2-\gamma}{4}\lambda_i\frac{\partial}{\partial\lambda_i}\int_{\mathbb{H}^1}\tilde{K}\left(\frac{\xi}{\lambda_i}+a_i\right)\frac{1}{\left((1+|x|^2+|y|^2)^2+t^2
\right)^2}\theta_0\wedge d\theta_0
\\&=-\frac{2-\gamma}{2}\Delta_{\theta_0} \tilde{K}(a_i)\int_{|\frac{\xi}{\lambda_i}|<\varepsilon}\frac{(x^2+y^2)\theta_0\wedge d\theta_0 }{\left((1+|x|^2+|y|^2)^2+t^2\right)^2} \\&\quad-\frac{2-\gamma}{4} \int_{|\frac{\xi}{\lambda_i}|\geqslant\varepsilon}\nabla_{\theta_0} \tilde{K}\left(\frac{\xi}{\lambda_i}+a_i\right)\frac{\xi}{\lambda_i}\frac{1}{\left((1+|x|^2+|y|^2)^2+t^2\right)^2}\theta_0\wedge d\theta_0 \\&\quad+o\left(\int_{|\frac{\xi}{\lambda_i}|<\varepsilon}|\frac{\xi}{\lambda_i}|^2
\frac{1}{\left((1+|x|^2+|y|^2)^2+t^2\right)^2}\right)\\&=-\frac{2-\gamma}{8}c_2\frac{\Delta_{\theta_0} \tilde{K}(a_i)}{\lambda_i^2}+o\left(\frac{1}{\lambda_i^2}\right).
\end{align*}\endd

\noindent{\textbf{Lemma B.3.2}}
For $i\neq j$, we have
\begin{align*}
&\int_{\mathbb{S}^3}K(\zeta) w_{g_j,\lambda_j}^
{\frac{2+\gamma}{2-\gamma}}\lambda_i\frac{\partial w_{g_i,\lambda_i}}{\partial\lambda_i}\theta_1
\wedge d\theta_1=O'(\lambda_i\frac{\partial\varepsilon_{ij}}{\partial\lambda_i})+O\left(\frac{1}{\lambda_i^
{2-\gamma}\lambda_j^{2+\gamma}}\right)+o\left(\frac{1}{\lambda_i^2}\right).
\end{align*}
\pf We have the following computations
\begin{align*}
&\int_{\mathbb{S}^3}K(\zeta) w_{g_j,\lambda_j}^
{\frac{2+\gamma}{2-\gamma}}\lambda_i\frac{\partial w_{g_i,\lambda_i}}{\partial\lambda_i}\theta_1
\wedge d\theta_1
\\&=\int_{\mathbb{H}^1}\tilde{K}(\xi)(\delta_{a_j,\lambda_j})^{\frac{2+\gamma}{2-\gamma}}\lambda_i\frac{\partial\delta_{a_i,\lambda_i}}{\partial\lambda_i}\theta_0\wedge d
\theta_0 \\&=\int_{\mathbb{H}^1}\tilde{K}(a_j)(\delta_{a_j,\lambda_j})^{\frac{2+\gamma}{2-\gamma}}\lambda_i\frac{\partial\delta_{a_i,\lambda_i}}{\partial\lambda_i}\theta_0\wedge d
\theta_0 \\&\quad+\int_{B(a_j,\varepsilon)}(\tilde{K}(\xi)-\tilde{K}(a_j))(\delta_{a_j,\lambda_j})^{\frac{2+\gamma}{2-\gamma}}\lambda_i\frac{\partial\delta_{a_i,\lambda_i}}{\partial
\lambda_i}\theta_0\wedge d\theta_0 \\&\quad+\int_{B^c(a_j,\varepsilon)}(\tilde{K}(\xi)-\tilde{K}(a_j))(\delta_{a_j,\lambda_j})^{\frac{2+\gamma}{2-\gamma}}\lambda_i\frac{\partial\delta_{a_i,\lambda_i}}{\partial
\lambda_i}\theta_0\wedge d\theta_0 \\&=O'(\lambda_i\frac{\partial\varepsilon_{ij}}{\partial\lambda_i})+o(\varepsilon_{ij})+O\left(\frac{1}{\lambda_i^
{2-\gamma}\lambda_j^{2+\gamma}}\right)+O(\varepsilon_{ij}^{\frac{2}{2-\gamma}}\ln\varepsilon_{ij}^{-1})+o\left(\frac{1}{\lambda_i^2}\right).
\end{align*}\endd

By some similar computations, we also get the following result.

\noindent{\textbf{Lemma B.3.3}}
For $i\neq j$, there holds
\begin{align*}
&\int_{\mathbb{S}^3}K(\zeta)( w_{g_j,\lambda_j})^{\frac{2\gamma}{2-\gamma}} w_{g_i,\lambda_i}\lambda_j\frac{\partial w_{g_j,\lambda_j}}{\partial\lambda_j}
\theta_1\wedge d\theta_1\\& =O'(\lambda_j\frac{\partial\varepsilon_{ij}}{\partial\lambda_j})+
O\left(\frac{1}{\lambda_i^{2-\gamma}\lambda_j^{2+\gamma}}\right)+O(\varepsilon_{ij}^{\frac{2}{2-\gamma}}\ln\varepsilon_{ij}^{-1})+o\left(\frac{1}{\lambda_i^2}\right).
\end{align*}

\noindent{\textbf{Lemma B.3.4}}
\begin{align*}
&\int_{\mathbb{S}^3}K(\zeta)\sum_{i\neq j,i\neq k}(\alpha_i w_{g_i,\lambda_i})^{\frac{2\gamma}{2-\gamma}}\textmd{inf}(\alpha_i w_{g_i,\lambda_i},
\alpha_k w_{g_k,\lambda_k})\lambda_j\frac{\partial w_{g_j,\lambda_j}}{\partial\lambda_j}\theta_1\wedge d\theta_1 \\&=
O\left(\sum_{i\neq j}\varepsilon_{ij}^{\frac{2}{2-\gamma}}\ln\varepsilon_{ij}^{-1}+\sum_{i\neq k}\varepsilon_{ik}^{\frac{2}{2-\gamma}}\ln\varepsilon_{ik}^{-1}\right).
\end{align*}
\pf From direct computation, we have $|\lambda_j\frac{\partial w_{g_j,\lambda_j}}{\partial\lambda_j}|\leqslant  w_{g_j,\lambda_j}.$ Then
\begin{align*}
&\left|\int_{\mathbb{S}^3}K(\zeta)\sum_{i\neq j,i\neq k}(\alpha_i w_{g_i,\lambda_i})^{\frac{2\gamma}{2-\gamma}}\textmd{inf}(\alpha_i w_{g_i,\lambda_i},
\alpha_k w_{g_k,\lambda_k})\lambda_j\frac{\partial w_{g_j,\lambda_j}}{\partial\lambda_j}\theta_1\wedge d\theta_1\right| \\&\leqslant \sum_{i\neq j,i\neq k}\int_{\mathbb{S}^3}(\alpha_i w_{g_i,\lambda_i})^{\frac{2\gamma}{2-\gamma}}(\alpha_i w_{g_i,\lambda_i})^{\frac{2-2\gamma}{2-\gamma}}(
\alpha_k w_{g_k,\lambda_k})^{\frac{\gamma}{2-\gamma}} w_{g_j,\lambda_j}\theta_1\wedge d\theta_1\\&\leqslant C\sum_{i\neq j,i\neq k}\int_{\mathbb{S}^3}(\alpha_i w_{g_i,\lambda_i})^{\frac{2}{2-\gamma}}\left((
\alpha_k w_{g_k,\lambda_k})^{\frac{2}{2-\gamma}}+ w_{g_j,\lambda_j}^{\frac{2}{2-\gamma}}\right)\theta_1\wedge d\theta_1\\&=
O\left(\sum_{i\neq j}\varepsilon_{ij}^{\frac{2}{2-\gamma}}\ln\varepsilon_{ij}^{-1}+\sum_{i\neq k}\varepsilon_{ik}^{\frac{2}{2-\gamma}}\ln\varepsilon_{ik}^{-1}\right).
\end{align*}\endd

Similarly, we have the following result, the proof is omitted.

\noindent{\textbf{Lemma B.3.5}}
\begin{align*}
&\int_{\mathbb{S}^3}K(\zeta)\sum_{i\neq j}(\alpha_j w_{g_j,\lambda_j})^{\frac{3\gamma-2}{2-\gamma}}\textmd{inf}((\alpha_i w_{g_i,\lambda_i})^2,
(\alpha_j w_{g_j,\lambda_j})^2)\lambda_j\frac{\partial w_{g_j,\lambda_j}}{\partial\lambda_j}\theta_1\wedge d\theta_1 \\&=
O\left(\sum_{i\neq j}\varepsilon_{ij}^{\frac{2}{2-\gamma}}\ln\varepsilon_{ij}^{-1}\right).
\end{align*}

By using the lemmas above, we have
\begin{align*}
J'(u)\left(\lambda_j\frac{\partial w_{g_j,\lambda_j}}{\partial\lambda_j}\right)&= \nonumber2\lambda(u)\left[-\sum_{i\neq j}O'(\lambda_j\frac{\partial\varepsilon_{ij}}{\partial\lambda_j})+
\frac{2-\gamma}{4}c_2\alpha_j\frac{\Delta_{\theta_1} K(g_j)}{ K(g_j)\lambda_j^2}\right.\\&\quad\left.o(\frac{1}{\lambda_j^2})+
o(\sum_{i\neq j} \varepsilon_{ij})\right].
\end{align*}
So we obtain the desired estimate of (\ref{es01}).

\subsection{Appendix C}

In this section, we take $W=\frac{1}{\lambda_j}\frac{\partial w_{g_j,\lambda_j}}{\partial g_j}$ in (\ref{J'}) and complete the proof of (\ref{es02}).
\begin{align*}
J'(u)\left(\frac{1}{\lambda_j}\frac{\partial w_{g_j,\lambda_j}}{\partial g_j}\right)&=
2\lambda(u)\left(\int_{\mathbb{S}^3}P_\gamma \left(\sum_{i=1}^p\alpha_i w_{g_i,\lambda_i}\right) \left(\frac{1}{\lambda_j}\frac{\partial w_{g_j,\lambda_j}}{\partial g_j}\right)\right.\\&\quad\left.-\lambda(u)^{\frac{2}{2-\gamma}}\int_{\mathbb{S}^3}K(\zeta)\left(\sum_{i=1}^p\alpha_i w_{g_i,
\lambda_i}\right)^
{\frac{2+\gamma}{2-\gamma}}\left(\frac{1}{\lambda_j}\frac{\partial w_{g_j,\lambda_j}}{\partial g_j}\right)\right).
\end{align*}
By the same reason of Lemma B.1, we have the following result.

\noindent{\textbf{Lemma C.1}} We have
\begin{align*}
\int_{\mathbb{S}^3}P_\gamma ( w_{g_j,\lambda_j})\frac{1}{\lambda_j}\frac{\partial w_{g_j,\lambda_j}}{\partial g_j}\theta_1\wedge d\theta_1 =0.
\end{align*}

\noindent{\textbf{Lemma C.2}} We have
\begin{align*}
\int_{\mathbb{S}^3}P_\gamma ( w_{g_i,\lambda_i})\frac{1}{\lambda_j}\frac{\partial w_{g_j,\lambda_j}}{\partial g_j}\theta_1\wedge d\theta_1
=O'(\frac{1}{\lambda_j}\frac{\partial\varepsilon_{ij}}{\partial g_j}).
\end{align*}
\Pf By some similar computations as in the proof of Lemma B.2, we can complete the proof of Lemma C.2. The details are omitted.\endd

\noindent{\textbf{Lemma C.3}} We have
\begin{align*}
&\int_{\mathbb{S}^3}K(\zeta)\left(\sum_{i=1}^p\alpha_i w_{g_i,\lambda_i}\right)^
{\frac{2+\gamma}{2-\gamma}}\left(\frac{1}{\lambda_j}\frac{\partial w_{g_j,\lambda_j}}{\partial g_j}\right)\\&=\alpha_j^
{\frac{2+\gamma}{2-\gamma}}\frac{2-\gamma}{2}c_2\frac{\nabla_{\theta_1}K(g_j)}{\lambda_j}+\sum_{i\neq j}O'(\frac{1}{\lambda_j}\frac{\partial\varepsilon_{ij}}{\partial g_j})+
O(\sum_{i\neq j} \varepsilon_{ij})+O\left(\frac{1}{\lambda_j^2}\right).
\end{align*}
\Pf Using (\ref{ineq5}), we have
\begin{align*}
&\left|\int_{\mathbb{S}^3}K(\zeta)\left(\sum_{i=1}^p\alpha_i w_{g_i,\lambda_i}\right)^
{\frac{2+\gamma}{2-\gamma}}\left(\frac{1}{\lambda_j}\frac{\partial w_{g_j,\lambda_j}}{\partial g_j}\right)\right.\\&\quad\left.-\int_{\mathbb{S}^3}K(\zeta)\left[
\sum_{i=1}^p(\alpha_i w_{g_i,\lambda_i})^
{\frac{2+\gamma}{2-\gamma}}+\frac{2+\gamma}{2-\gamma}(\alpha_j w_{g_j,\lambda_j})^
{\frac{2\gamma}{2-\gamma}}(\sum_{i\neq j}\alpha_i w_{g_i,\lambda_i})\right]\frac{1}{\lambda_j}\frac{\partial w_{g_j,\lambda_j}}{\partial g_j}\right|\\&\leqslant M_5\left(\int_{\mathbb{S}^3}K(\zeta)\sum_{i\neq j,i\neq k}(\alpha_i w_{g_i,\lambda_i})^{\frac{2\gamma}{2-\gamma}}\textmd{inf}(\alpha_i w_{g_i,\lambda_i},
\alpha_k w_{g_k,\lambda_k})\frac{1}{\lambda_j}\frac{\partial w_{g_j,\lambda_j}}{\partial g_j}\right)\\&\quad+M_5 \left(\int_{\mathbb{S}^3}K(\zeta)\sum_{i\neq j}(\alpha_j w_{g_j,\lambda_j})^
{\frac{3\gamma-2}{2-\gamma}}\textmd{inf}((\alpha_i w_{g_i,\lambda_i})^2,
(\alpha_j w_{g_j,\lambda_j})^2)\frac{1}{\lambda_j}\frac{\partial w_{g_j,\lambda_j}}{\partial g_j}\right).
\end{align*}

Together with Lemma C.3.1, Lemma C.3.2, Lemma C.3.3, Lemma C.3.4, and Lemma C.3.5, we get the desired results.\endd

\noindent{\textbf{Lemma C.3.1}} We have
\begin{align*}
\int_{\mathbb{S}^3}K(\zeta)(w_{g_j,\lambda_j})^{\frac{2+\gamma}{2-\gamma}}\frac{1}{\lambda_j}\frac{\partial w_{g_j,\lambda_j}}{\partial g_j}\theta_1\wedge d\theta_1 =
\frac{2-\gamma}{4}c_2\frac{\nabla_{\theta_1}K(g_j)}{\lambda_j}+O(\frac{1}{\lambda_j^2}).
\end{align*}
\Pf We understand the vectors in the following formula as row vectors, then there holds
\begin{align*}
&\int_{\mathbb{S}^3}K(\zeta)(w_{g_j,\lambda_j})^{\frac{2+\gamma}{2-\gamma}}\frac{1}{\lambda_j}\frac{\partial w_{g_j,\lambda_j}}{\partial g_j}\theta_1\wedge d\theta_1\\
&=\left[\int_{\mathbb{H}^1}\tilde{K}(\xi)(\delta_{a_j,\lambda_j})^{\frac{2+\gamma}{2-\gamma}}\frac{1}{\lambda_j}\frac{\partial\delta_{a_j,\lambda_j}}{\partial a_j}\theta_0\wedge d\theta_0\right]d \mathcal{C}^{-1}|_{g_j}.
\end{align*}

On the right hand side, we have
\begin{align*}
&\int_{\mathbb{H}^1}\tilde{K}(\xi)(\delta_{a_j,\lambda_j})^{\frac{2+\gamma}{2-\gamma}}\frac{1}{\lambda_j}\frac{\partial\delta_{a_j,\lambda_j}}{\partial a_j}\theta_0\wedge d\theta_0 \\&=\int_{B(a_j,\varepsilon)}(\tilde{K}(\xi)-\tilde{K}(a_j))(\delta_{a_j,\lambda_j})^{\frac{2+\gamma}{2-\gamma}}\frac{1}{\lambda_j}
\frac{\partial\delta_{a_j,\lambda_j}}{\partial a_j}\\&\quad+\int_{B^c(a_j,\varepsilon)}(\tilde{K}(\xi)-\tilde{K}(a_j))(\delta_{a_j,\lambda_j})^
{\frac{2+\gamma}{2-\gamma}}\frac{1}{\lambda_j}\frac{\partial\delta_{a_j,\lambda_j}}{\partial a_j}\\&=
\int_{B(a_j,\varepsilon)}\nabla_{\theta_0}\tilde{K}(a_j)(\xi-a_j)(\delta_{a_j,\lambda_j})^{\frac{2+\gamma}{2-\gamma}}\frac{1}{\lambda_j}
\frac{\partial\delta_{a_j,\lambda_j}}{\partial a_j}+O\left(\int_{B(a_j,\varepsilon)}(\delta_{a_j,\lambda_j})^{\frac{5}{2-\gamma}}|\xi-a_j|^3\right)\\&\quad+\int_{B^c(a_j,\varepsilon)}(\tilde{K}(\xi)-\tilde{K}(a_j))(\delta_{a_j,\lambda_j})^
{\frac{2+\gamma}{2-\gamma}}\frac{1}{\lambda_j}\frac{\partial\delta_{a_j,\lambda_j}}{\partial a_j}\theta_0\wedge(d\theta_0)^n\\&=
\frac{2-\gamma}{4}c_2\frac{\nabla_{\theta_0}\tilde{K}(a_j)}{\lambda_j}+O(\frac{1}{\lambda_j^2}).
\end{align*}\endd

\noindent{\textbf{Lemma C.3.2}}
For $i\neq j$, we have
\begin{align*}
\int_{\mathbb{S}^3}K(\zeta)(w_{g_i,\lambda_i})^{\frac{2+\gamma}{2-\gamma}}\frac{1}{\lambda_j}\frac{\partial w_{g_j,\lambda_j}}{\partial g_j}\theta_1\wedge d\theta_1 =O'(\frac{1}{\lambda_j}\frac{\partial\varepsilon_{ij}}{\partial g_j})+o(\varepsilon_{ij})+O(\frac{1}{\lambda_i^
{2+\gamma}\lambda_j^{2-\gamma}}).
\end{align*}
\Pf Similarly as in the proof of Lemma C.3.1, we have
\begin{align*}
&\int_{\mathbb{S}^3}K(\zeta)(w_{g_i,\lambda_i})^{\frac{2+\gamma}{2-\gamma}}\frac{1}{\lambda_j}\frac{\partial w_{g_j,\lambda_j}}{\partial g_j}\theta_1\wedge d\theta_1\\&=\left[\int_{\mathbb{H}^1}\tilde{K}(\xi)(\delta_{a_i,\lambda_i})^{\frac{2+\gamma}{2-\gamma}}\frac{1}{\lambda_j}\frac{\partial\delta_{a_j,\lambda_j}}{\partial a_j}\theta_0\wedge d\theta_0\right]d \mathcal{C}^{-1}|_{g_j},
\end{align*}and
\begin{align*}
&\int_{\mathbb{H}^1}\tilde{K}(\xi)(\delta_{a_i,\lambda_i})^{\frac{2+\gamma}{2-\gamma}}\frac{1}{\lambda_j}\frac{\partial\delta_{a_j,\lambda_j}}{\partial a_j}\theta_0\wedge d\theta_0 \\&=\tilde{K}(a_i)\frac{1}{\lambda_j}\frac{\partial}{\partial a_j}\int_{\mathbb{H}^1}(\delta_{a_i,\lambda_i})^{\frac{2+\gamma}{2-\gamma}}\delta_{a_j,\lambda_j}\theta_0\wedge d
\theta_0 \\&\quad+\int_{B(a_i,\varepsilon)}(\tilde{K}(\xi)-\tilde{K}(a_i))(\delta_{a_i,\lambda_i})^{\frac{2+\gamma}{2-\gamma}}\frac{1}{\lambda_j}\frac{\partial\delta_{a_j,\lambda_j}}{\partial a_j}\theta_0\wedge d\theta_0 \\&\quad+\int_{B^c(a_i,\varepsilon)}(\tilde{K}(\xi)-\tilde{K}(a_i))(\delta_{a_i,\lambda_i})^{\frac{2+\gamma}{2-\gamma}}\frac{1}{\lambda_j}\frac{\partial\delta_{a_j,\lambda_j}}{\partial a_j}\theta_0\wedge d\theta_0 \\&=O'(\frac{1}{\lambda_j}\frac{\partial\varepsilon_{ij}}{\partial a_j})+o(\varepsilon_{ij})+O\left(\frac{1}{\lambda_i^
{2+\gamma}\lambda_j^{2-\gamma}}\right).
\end{align*}\endd

Similarly we have the following result, but its proof is more simple.

\noindent{\textbf{Lemma C.3.3}}
For $i\neq j$, we have
\begin{align*}
&\int_{\mathbb{S}^3}K(\zeta)( w_{g_j,\lambda_j})^{\frac{2\gamma}{2-\gamma}} w_{g_i,\lambda_i}\frac{1}{\lambda_j}\frac{\partial w_{g_j,\lambda_j}}{\partial g_j}\theta_0\wedge d\theta_0 =O(\varepsilon_{ij}).
\end{align*}

\noindent{\textbf{Lemma C.3.4}} We have
\begin{align*}
&\int_{\mathbb{S}^3}K(\zeta)\sum_{i\neq j,i\neq k}(\alpha_i w_{g_i,\lambda_i})^{\frac{2\gamma}{2-\gamma}}\textmd{inf}(\alpha_i w_{g_i,\lambda_i},
\alpha_k w_{g_k,\lambda_k})\frac{1}{\lambda_j}\frac{\partial w_{g_j,\lambda_j}}{\partial g_j}\theta_0\wedge d\theta_0 \\&=
O\left(\sum_{i\neq j}\varepsilon_{ij}^{\frac{2}{2-\gamma}}\ln\varepsilon_{ij}^{-1}+\sum_{i\neq k}\varepsilon_{ik}^{\frac{2}{2-\gamma}}\ln\varepsilon_{ik}^{-1}\right).
\end{align*}
\pf Through direct computation, we have $|\frac{1}{\lambda_j}\frac{\partial w_{g_j,\lambda_j}}{\partial g_j}|\leqslant  w_{g_j,\lambda_j}.$ Using Lemma A.2, we get the result.\endd

We also have the following result, its proof is omitted.

\noindent{\textbf{Lemma C.3.5}} We have
\begin{align*}
&\int_{\mathbb{S}^3}K(\zeta)\sum_{i\neq j}(\alpha_j w_{g_j,\lambda_j})^{\frac{3\gamma-2}{2-\gamma}}\textmd{inf}((\alpha_i w_{g_i,\lambda_i})^2,
(\alpha_j w_{g_j,\lambda_j})^2)\frac{1}{\lambda_j}\frac{\partial w_{g_j,\lambda_j}}{\partial g_j}\theta_0\wedge d\theta_0 \\&=
O\left(\sum_{i\neq j}\varepsilon_{ij}^{\frac{2}{2-\gamma}}\ln\varepsilon_{ij}^{-1}\right).
\end{align*}

By using the lemmas above, we have
\begin{align*}
J'(u)\left(\frac{1}{\lambda_j}\frac{\partial w_{g_j,\lambda_j}}{\partial g_j}\right)=-2\lambda(u)\alpha_j
c_2\frac{\nabla_{\theta_1}K(g_j)}{ K(g_j)\lambda_j}+
O(\sum_{i\neq j} \varepsilon_{ij}+\frac{1}{\lambda_j^2}).
\end{align*}
This is the desired estimate of (\ref{es02}).


\end{document}